\documentclass{amsart}
\usepackage[utf8]{inputenc}

\usepackage{geometry}
\usepackage{epsdice}

\usepackage{enumitem}

\usepackage{color}

\usepackage[OT2,T1]{fontenc}

\usepackage[final]{microtype}

\usepackage{tikz}
\usepackage{tikz-cd}
\usepackage{xypic}

\usepackage[draft=false,linktocpage=true]{hyperref}
\usepackage[capitalise]{cleveref}

\usepackage{relsize}
\usepackage[bbgreekl]{mathbbol}
\usepackage{amsfonts}
\DeclareSymbolFontAlphabet{\mathbb}{AMSb} 
\DeclareSymbolFontAlphabet{\mathbbl}{bbold}

\usepackage{mathtools}
\usepackage{amsmath}
\usepackage{mathrsfs}

\usepackage{amsthm,amssymb}
\numberwithin{equation}{section}

\theoremstyle{plain}
\newtheorem{theorem}[equation]{Theorem}

\newtheorem{proposition}[equation]{Proposition}
\newtheorem{lemma}[equation]{Lemma}

\newtheorem{corollary}[equation]{Corollary}

\theoremstyle{definition}
\newtheorem{definition}[equation]{Definition}
\newtheorem{notation}[equation]{Notation}
\newtheorem{construction}[equation]{Construction}

\newtheorem{example}[equation]{Example}

\newtheorem{remark}[equation]{Remark}

\setlist[enumerate]{leftmargin=*}
\setlist[itemize]{leftmargin=*}

\newcommand{\conj}{\mathrm{conj}}
\newcommand{\Fil}{\mathrm{Fil}}
\newcommand{\qsyn}{\mathrm{qSyn}}

\newcommand{\w}{\text}

\date{First version August 2, 2021; Revised version April 20, 2022.}

\begin{document}

\title{On Endomorphisms of the de Rham cohomology functor}

\author{Shizhang Li}

\address{Morningside Center of Mathematics and Hua Loo-Keng Key Laboratory of Mathematics,
Academy of Mathematics and Systems Science, Chinese Academy of Sciences, 
Beijing 100190, China}

\email{lishizhang@amss.ac.cn}

\author{Shubhodip Mondal}

\address{Department of Mathematics, University of Michigan, 530 Church Street,
  Ann Arbor, MI 48109}

\email{smondal@umich.edu}


\maketitle
\begin{abstract}
We compute the moduli of endomorphisms of the de Rham and crystalline cohomology functors, viewed
as a cohomology theory on smooth schemes over truncated Witt vectors. 
As applications of our result, we deduce Drinfeld's refinement of the classical Deligne--Illusie decomposition result for de Rham cohomology of varieties in characteristic $p>0$ that are liftable to $W_2$, 
and prove further functorial improvements.
\end{abstract}

\tableofcontents
\newpage
\section{Introduction}
Let $A$ be a ring and let $X$ be a smooth $A$-scheme. 
The algebraic de Rham cohomology is a cohomology theory designed by Grothendieck. 
It is defined functorially by sending $X$ to the hypercohomology of the de Rham complex $\Omega_{X/A}^*.$ 
The de Rham complex $\Omega_{X/A}^*$ is not just a complex, but also has an additional structure of a sheaf of commutative differential graded algebras. 
One can therefore view the output of de Rham cohomology as a commutative algebra object in the derived $\infty$-category $D(A)$
which we denote by $\w{CAlg}(D(A)).$ 
This way, one obtains a functor $\w{dR}_{(\cdot)/A}: \w{Alg}^{\w{sm}}_A \to \w{CAlg}(D(A)),$ which sends any smooth $A$-algebra $R \mapsto \mathrm{dR}_{R/A} \in \w{CAlg}(D(A)).$ The primary goal of our paper is to study endomorphisms of this functor.

Studying properties of the de Rham cohomology theory as a functor is interesting for a number of reasons. From a technical point of view, in certain situations, showing that the de Rham cohomology functor has no nontrivial automorphisms has been used as a key tool in \cite{BLM20} and \cite{LL20} to prove that certain constructions are functorially isomorphic. 
Further, in \cite{Mon21}, 
it was shown that one can reconstruct the theory of crystalline cohomology as the 
\emph{unique} deformation of de Rham cohomology theory viewed as a functor defined on smooth $\mathbb{F}_p$-schemes. 

From a different perspective, any property enjoyed by the de Rham cohomology functor will in particular be enjoyed by de Rham cohomology of every smooth algebraic variety.
For example, if the functor $\w{dR}_{(\cdot)/A}$
has many endomorphisms, one potentially obtains many interesting endomorphisms of the de Rham cohomology of any smooth algebraic variety, 
which could be useful for making interesting geometric conclusions. 
The classical study and usage of Frobenius operator on de Rham or crystalline cohomology theory is 
an instance of such a perspective.


In this paper, our main motivating questions are the following, which can be seen as a ``moduli'' enhancement of the question of endomorphisms of the de Rham cohomology functor:
\begin{enumerate}
    \item Given a ring $A$, consider the functor $\mathrm{dR}$ that sends any smooth $A$-algebra $R \mapsto \mathrm{dR}_{R/A} \in \w{CAlg}(D(A))$, 
    what is the endomorphism monoid\footnote{\textit{A priori} we get a monoid object in spaces rather than an actual monoid.
But in the cases of interest to us, this space is discrete, see \Cref{discrete endo space} and \Cref{dRmn discrete endo space}.} of this functor? \vspace{1mm}

    \item More generally, let $B$ be an arbitrary $A$-algebra, consider the analogous functor $R \mapsto \mathrm{dR}_{R/A} \otimes_{A} B \in \w{CAlg}(D(B))$,
    what is the endomorphism monoid of this functor? \vspace{1mm}
    \item Finally consider the presheaf\footnote{Mathew pointed out to us that this presheaf is automatically
    an fpqc sheaf by flat descent.} (of monoids) on $\mathrm{Alg}_A^{\w{op}}$ sending an $A$-algebra $B$ to the endomorphism
    monoid in previous question. Is it represented by a (monoid) scheme?
    If so, what is the representing monoid scheme?
\end{enumerate}




In this paper, we address the above questions when $A= W_n(k)$ for any perfect ring $k$, 
where $W_n(k)$ denotes the ring of $n$-truncated Witt vectors. We expect the methods to be extendable to more general base rings 
but we do not pursue that direction further in this paper. 

\subsection*{A foretaste of the Main Theorem}
For simplicity, for now let us focus on the case where $A = \mathbb{Z}/p^n$ or
$\mathbb{Z}_p$ and $B$ is an $\mathbb{F}_p$-algebra.

\begin{theorem}[Special case of Main \Cref{representing endomorphism monoid}]
\label{mainintro}
\leavevmode
\begin{enumerate}
\item When $A = \mathbb{F}_p$, the endomorphism monoid of $\mathrm{dR}_{(\cdot)/A} \otimes_A B$ is $\underline{\mathbb{N}}(\mathrm{Spec}(B))$
where $\underline{\mathbb{N}}$ denotes the constant monoid scheme associated with the natural numbers.
\item However, when $A = \mathbb{Z}/p^n$ for $n \geq 2$, the endomorphism monoid of $\mathrm{dR}_{(\cdot)/A} \otimes_A B$ is a semi-direct product of $\underline{\mathbb{N}}(\mathrm{Spec}(B))$ with a group $W(B)^{\times}[F]$, where the latter denotes the Frobenius kernel
of the unit group in $W(B)$.
\end{enumerate}
\end{theorem}

\begin{remark}
\leavevmode
\begin{enumerate}
\item Roughly speaking, when $A = \mathbb{Z}/p^n$ for $n \geq 2$, \cref{mainintro} says that the endomorphism monoid of
$\mathrm{dR}$ is very large. 
More precisely, \cref{mainintro} in particular provides an action of $W^\times [F]$ on the mod $p$ de Rham cohomology of a variety liftable to $W_2.$
Recently, Drinfeld has also observed an action of $W^{\times}[F]$ on the mod $p$ de Rham cohomology, using the theory of ``prismatization'' due to him and independently due to Bhatt--Lurie. The main new ingredient of \cref{mainintro} is to go \emph{beyond} this action and \emph{classify all the endomorphisms}. Interestingly, our proof of \cref{mainintro} does not make any use of prismatization and only uses the stacky approach to de Rham cohomology theory in positive characteristic that already appeared in \cite{drinfeld2018stacky}. However, while the stacky approach (including the theory of prismatization) helps in constructing the endomorphisms, it does not \textit{apriori} offer any strategy to prove that they are \textit{all} the endomorphisms. To achieve this, we employ some very different additional techniques in the proof of \cref{mainintro}, such as the theory of affine stacks due to To\"en \cite{MR2244263}, a version of the topos theoretic cotangent complex (see \cref{appendix}) due to Illusie \cite{Ill} and some explicit computations when necessary.

\item The $W^\times [F]$ action resulting from \cref{mainintro} will be utilized to prove a strengthened version
of the Deligne--Illusie decomposition, see \Cref{corintro}. See \cref{nofunctsplit} for an application of the full classification offered by \cref{mainintro}.
\item From the above calculation one finds that for $A = \mathbb{Z}/p^n$, the association $B \mapsto \mathrm{End}(\mathrm{dR}_{(\cdot)/A} \otimes_A B)$
defines a sheaf of monoids representable by a scheme denoted by $\mathrm{End}_{1,n}$.
The representing monoid scheme depends on $A= \mathbb{Z}/p^n$ and stabilizes
when $n \ge 2$. 
\end{enumerate}
\end{remark}

The stabilization we refer to means the following. Note that we have a natural commutative diagram
\[
\xymatrix{
\text{Alg}_{\mathbb{Z}/p^n}^{\mathrm{sm}} \ar[rd]_{\mathrm{dR} \otimes B} \ar[rr]^{\text{mod } p^{n-1}} & & \text{Alg}_{\mathbb{Z}/p^{n-1}}^{\mathrm{sm}} \ar[ld]^{\mathrm{dR} \otimes B} \\
& \w{CAlg}(D(B)), &
}
\]
which induces a sequence of maps of schemes
\[
\mathrm{End}_{1,1} \to 
\mathrm{End}_{1,2}  \to
\ldots \to \mathrm{End}_{1,n} \to \ldots.
\]
Our theorem says that the first map is a closed immersion and all subsequent maps are isomorphisms.

\begin{remark}
We see that the representing monoid scheme stabilizes as soon
as $A$ leaves characteristic $p$, this indicates the functorial Frobenius endomorphism
is solely responsible for the rigidity of de Rham
cohomology theory in characteristic $p$.
\end{remark}

Regarding endomorphisms of de Rham cohomology itself, we also get the following:
\begin{theorem}[{Special case of \Cref{integral endomorphism prop}}]
\label{p-adic dR has no auto}
When $A = \mathbb{Z}_p$, 
the endomorphism monoid of $\mathrm{dR}_{(\cdot)/A}^{\wedge}$ is $\mathbb{N}$, given by powers of the Frobenius.
\end{theorem}

Here the $\mathrm{dR}^{\wedge}$ denotes the $p$-adic derived de Rham cohomology theory,
\textit{c.f.}~\cite{Bha12}.
The fact that there is no automorphism of $p$-adic derived de Rham cohomology theory when the base ring
is $p$-complete and $p$-torsionfree was observed in \cite[Theorem 3.13]{LL20}.

\begin{remark}
In both cases of $A = \mathbb{F}_p$ and $\mathbb{Z}_p$ above, we only see powers of the Frobenius
as endomorphisms of the ($p$-adic) de Rham cohomology.
But they are due to two different reasons:
when $A = \mathbb{F}_p$ it is due to existence of the Frobenius endomorphism on the category
of $A$-algebras, whereas for $A = \mathbb{Z}_p$ it comes from the fact that
$A$ is $p$-torsionfree, so a certain huge group scheme has no nontrivial $A$-valued point in this case.
\end{remark}

In \Cref{representing endomorphism monoid} we work in a more general setting, 
namely we calculate the moduli of endomorphisms of crystalline cohomology theory, leading to the sheaves $\mathrm{End}_{m,n}$
(see \cref{corend} for the precise statement).
The result is similar:
the Frobenius endomorphisms that people ``knew and loved'' corresponds to the monoid underlying the connected components of the whole endomorphism monoid; in fact there is a distinguished point in each component which corresponds to a power of the Frobenius endomorphism. 
Furthermore, the identity component also stabilizes to a large and mysterious group scheme (see \cref{notation of aut gp})
which demands further investigations (see \cref{twists}).
One surprising feature is that the above group scheme is 
non-flat over the base in the general setting of crystalline cohomology.

\subsection*{Application to the Deligne--Illusie decomposition}
As an application of the $W^\times [F]$-action, Drinfeld observed a refinement of the Deligne--Illusie decomposition, 
which was communicated to us by Bhatt
(see \cite[Example 4.7.17, Remark 4.7.18]{BhaLur} and \cite[Remark 5.16]{BhaLur2}):
since $\mu_p \subset W^{\times}[F]$, the mod $p$ de Rham cohomology of varieties liftable to $W_2$ has the structure of a $\mu_p$-representation.
It is easy to see that the $W^{\times}[F]$-action preserves conjugate filtration. Then one needs to show that the $i$-th graded piece
of conjugate filtration is pure of weight $i \in \mathbb{Z}/p$ as a $\mu_p$-representation. 
In \cite{BhaLur} and \cite{BhaLur2}, 
this statement is proven by establishing a relation between $W^\times [F]$-action and the ``Sen operator''
defined in loc.~cit.
In \cref{computing mup action}, we use a more direct argument to check that the weight statement holds for the $W^\times[F]$-action coming from our \cref{mainintro}.


\cref{mainintro} coupled with the calculation of weights from \cref{computing mup action} as above immediately implies the following improvement of results due to Achinger and Suh \cite[Theorem 1.1]{AchSuh}, which in turn
is a strengthening of Deligne--Illusie's result \cite[Corollaire 2.4]{DI87}. 
In particular, our approach gives a proof, different from Bhatt--Lurie's, of the following result which does not make any use of prismatization.
\begin{theorem}[Drinfeld, \cite{BhaLur2}, see \Cref{Drinfeld splitting}]
\label{corintro}
Let $k$ be a perfect ring of characteristic $p > 0$, let $X$ be a smooth scheme over $W_2(k)$,
and let $a \leq b \leq a + p-1$.
Then the canonical truncation $\tau_{[a,b]}(\Omega^{\bullet}_{X_k/k})$ splits.
Moreover the splitting is functorial in the lift $X$ of $X_k$.
\end{theorem}

Since our calculation shows that the endomorphism monoid of mod $p$ de Rham cohomology stabilizes after $W_2$,
philosophically it says that further liftability over $W_n$ for $n > 2$ provides no extra knowledge on the mod $p$ de Rham cohomology.

It is still an open problem whether there exists a smooth variety $X$ (necessarily
of dimension $\dim X > p$) over $k$ which lifts to $W_2(k)$, for which the de Rham complex is not decomposable. Using \cref{mainintro}, we obtain a somewhat negative result in this direction: we show that de Rham complex of smooth varieties over $k$ liftable to $W_2(k)$ 
does not completely decompose in a \emph{functorial} manner as a commutative algebra object in the derived category. 

\begin{corollary}[see \Cref{no further functorial splitting}]\label{nofunctsplit}
There is no functorial splitting
\[
\mathrm{dR}_{(- \otimes_{W_2(k)} k)/k} \simeq
\bigoplus_{i \in \mathbb{N}_{\geq 0}} \mathrm{Gr}^{\mathrm{conj}}_{i}\big(\mathrm{dR}_{(- \otimes_{W_2(k)} k)/k}\big)
\]
as a functor from smooth $W_2(k)$-algebras to $\mathrm{CAlg}(D(k))$.
\end{corollary}
The above statement was also observed by Mathew. 
His idea for a proof does not use the full calculation of endomorphism monoid as in \cref{mainintro},
whereas for us it is a consequence of that calculation.

Lastly one may wonder if the Drinfeld splitting agrees with the Deligne--Illusie splitting (which has an $\infty$-categorical functorial
enhancement~\cite[Theorem 1.3.21 and Proposition 1.3.22]{KP19}).
Both splittings are obtained from the splitting of the first conjugate filtration,
via an averaging process, \textit{c.f.}~the step (a) in proof of \cite[Th\'{e}or\`{e}me 2.1]{DI87}.
To guarantee the above two splittings are functorially the same, we show the following uniqueness. 

\begin{theorem}[see \Cref{uniqueness of functorial splitting of Conj1} for the precise statement]
There is a unique functorial splitting (as $\Fil^{\mathrm{conj}}_{0}$-modules)
\[
\Fil^{\mathrm{conj}}_{1}(\mathrm{dR}_{(- \otimes_{W_2(k)} k)/k}) = 
\Fil^{\mathrm{conj}}_{0}(\mathrm{dR}_{(- \otimes_{W_2(k)} k)/k}) \oplus \mathrm{Gr}^{\mathrm{conj}}_{1}(\mathrm{dR}_{(- \otimes_{W_2(k)} k)/k}).
\]
\end{theorem}

In particular the Deligne--Illusie splitting in \cite{KP19}, the Drinfeld splitting in 
\cite{BhaLur} \cite{BhaLur2} as well as the splitting induced by \cref{mainintro}
must all agree.
\vspace{2mm}
\noindent
\subsection*{Outline of the proof of Theorem 1.1}\label{outlineproof} Let us briefly outline the key ingredients in the proof of \cref{mainintro}. 
In doing so, we will also give a rough outline of the paper.
For simplicity, let us fix $A = \mathbb{Z}/p^n,$ and let $B$ be an $\mathbb{F}_p$-algebra.
\begin{enumerate}
\item[(0)] The statements \cref{mainintro}.(1) as well as \Cref{p-adic dR has no auto} is within reach of the quasi-syntomic descent techniques
introduced in \cite{BMS2}, see \Cref{section 3}.
We also use quasi-syntomic descent techniques to show that the endomorphism spaces of interest to us are actually discrete
(see \cref{discrete endo space}, \cref{dRmn discrete endo space}).\vspace{1mm}
\item But for \cref{mainintro}.(2), we need to make use of the stacky approach to de Rham or crystalline cohomology due to Drinfeld \cite{drinfeld2018stacky}, \cite{prismatization}, which can be seen as a positive characteristic analogue of Simpson's de Rham stack \cite{Si96}. For our paper, we use a compressed version of the stacky approach: the functor $\w{dR}_{(\cdot)/A} \otimes _A B$ is built as  \textit{unwinding} (see \cref{unwindringstack}) of an $A$-algebra stack over $B;$ this stack is denoted by $\mathbb{A}^{1, \w{dR}}_B$ (we often omit $B$ to ease the notation). This unwinding construction is a variant of a construction used in \cite[\S 3]{Mon21}.
Note the amusing switch of roles played by $A$ and $B$: the de Rham cohomology theory
is a cohomology theory for varieties over $A$ with coefficient ring being $B$,
whereas the stack $\mathbb{A}^{1, \w{dR}}_B$ is an $A$-algebra object over $B$.\vspace{1mm}

\item It turns out that the underlying stack $\mathbb{A}^{1,\w{dR}}$ is an affine stack, in the sense of To\"{e}n \cite{MR2244263}. Roughly speaking, for affine stacks one can pass to the ``ring'' of derived global sections in a lossless manner. Using this property, in \cref{unwindequiv}, we show that $\w{End}(\w{dR}_{(\cdot)/A} \otimes _A B) \simeq \w{End}_{A\w{-Alg-St}}(\mathbb{A}^{1, \w{dR}}_B).$
Here the latter endomorphism is taken in the category of $A$-algebra stacks over $B.$ \vspace{1mm}

\item Using the description of $\mathbb{A}^{1,\w{dR}}$ as the quotient stack $[W/pW],$ where $W$ denotes the ring scheme of $p$-typical Witt vectors, in \cref{enough} we construct ``enough'' endomorphisms of $\w{dR}_{(\cdot)/A} \otimes_A B$ and show that the endomorphism monoid is at least as big as \cref{mainintro} claims. \vspace{1mm}

\item To finish the proof of \cref{mainintro}, one needs to show that there is no more endomorphism than the ones already constructed. To do so, we interpret an endomorphism of the algebra stack $\mathbb{A}^{1,\w{dR}}$ as a deformation of an endomorphism of the sheaf of rings $\pi_0 (\mathbb{A}^{1,\w{dR}}).$ We know that (see \cref{cohomology of A1dR}) $\pi_0 (\mathbb{A}^{1,\w{dR}}_B)=\mathbb{G}_{a,B}$ because $B$ is an $\mathbb{F}_p$-algebra.
Then we use the formalism of topos theoretic cotangent complex due to Illusie \cite{Ill} (see \Cref{appendix}) to understand this deformation problem. This is carried out in \cref{representing endomorphism monoid}, where we use the cotangent complex and the transitivity triangle to finish calculating
the desired endomorphism monoid.
\end{enumerate}{}

\begin{remark}
Let $A$ and $B$ be as above. Combining the steps mentioned in (2) and (3) above, one obtains that an endomorphism of the functor $\w{End}(\w{dR}_{(\cdot)/A} \otimes _A B)$ is the same datum as a natural endomorphism of $W(S)/^{L} p,$ as an animated $A$-algebra, for every (discrete) $B$-algebra $S$.
\end{remark}{}

\subsection*{Acknowledgements} We are grateful to Bhargav Bhatt for informing us about \cref{corintro} as well as many 
stimulating discussions and helpful suggestions,
and for organizing an informal seminar during the summer of 2021 in which we presented this work to the participants: Attilio Castano, Haoyang Guo, Andy Jiang, Emanuel Reinecke, Gleb Terentiuk, Jakub Witaszek,
and Bogdan Zavyalov; we thank them for their interest and helpful conversations. We are also very thankful to Piotr Achinger, Ben Antieau, Johan de Jong, Luc Illusie, Dmitry Kubrak and Akhil Mathew for their comments and valuable feedbacks. Special thanks to the anonymous referee for many comments and suggestions on the paper.

The first named author thanks the support of AMS-Simons Travel Grant 2021-2023. 
The second named author thanks the support of NSF grant DMS \#1801689 through Bhargav Bhatt.

\section{Stacky approach to de Rham cohomology}
The goal of this section is to describe the stacky approach to de Rham cohomology theory due to Drinfeld \cite{drinfeld2018stacky}. Roughly, given a scheme $X,$ Drinfeld constructed a stack $X^{\w{dR}}$ such that $R\Gamma(X^{\w{dR}}, \mathcal{O})$ recovers the de Rham cohomology $R\Gamma_{\w{dR}}(X).$ This should be seen as a positive characteristic variant of the earlier construction of the de Rham stack due to Simpson \cite{Si96}.

For our purposes, we will need to work with a certain compressed version of this construction. Our goal is to consider a \textit{single} stack with enough structure encoded, which can naturally ``unwind'' itself to construct the stack $X^{\w{dR}}$ for \textit{every} scheme $X.$ To this end, we will begin by discussing quasi-ideals (\textit{c.f.}~\cite[\S 3.1]{prismatization}, ~\cite[\S 3.2]{Mon21}) and ring stacks, which formulates exactly the kind of extra structures on a stack one needs to work with in order to use the unwinding machine. After that, we will discuss the construction of this unwinding functor and explain how to build a cohomology theory from a ring stack in general. We will then discuss the particular ring stack $\mathbb{A}^{1, \w{dR}}$ which gives rise to de Rham cohomology theory via this construction. For later application, the fact that the stack $\mathbb{A}^{1,\w{dR}}$ is an \textit{affine stack} in the sense of \cite{MR2244263} will be of particular importance to us. Therefore, we will record the relevant definitions in this section as well.

\subsection{Quasi-ideals}\label{sec2.1}


\begin{definition}[Quasi-ideals]\label{qid}Let $R$ be a ring and $M$ be an $R$-module equipped with a map $d: M \to R$ of $R$-modules which satisfies $d(x) \cdot y = d(y) \cdot x$ for any pair $x,y \in M$.
Such a data $d: M \to R$ satisfying the aforementioned condition will be called a \textit{quasi-ideal} in $R$ or simply a quasi-ideal. 

A morphism of quasi-ideals $(d_1: M_1 \to R_1)\to (d_2: M_2 \to R_2)$ is defined to be a pair of maps $a: M_1 \to M_2$ and $b: R_1 \to R_2$ such that the following compatibilities hold:

\begin{enumerate}
\item $d_2 a = b d_1.$
\item $a (r_1 m_1) = b(r_1) a(m_1).$
\item $b$ is a ring homomorphism.
\item $a$ is linear.

\end{enumerate}
\noindent

In other words, we want a commutative diagram as below: \begin{center}
    \begin{tikzcd}
M_1 \arrow[d, "d_1"] \arrow[r, "a"] & M_2 \arrow[d, "d_2"] \\
R_1 \arrow[r, "b"]                  & R_2                 
\end{tikzcd}
\end{center} such that $b$ is a ring homomorphism and $a$ is an $R_1$-module map $M_1 \to b_* M_2.$ The category of quasi-ideals will be denoted as $\w{QID}.$ 



\begin{construction}[Quasi-ideal as a simplicial abelian group]Given a quasi-ideal $(d \colon M \to R ),$ we obtain a map 
$t \colon T \coloneqq M \times R \to R$ given by $(m,r) \to r+ d(m).$ There is another map $s \colon M \times R \to R$ given by $(m,r) \to r.$ There is also a degeneracy map $e \colon R \to M \times R$ given by $r \to (0,r).$ Lastly, there is a map $c \colon T \times_{R, s, t} T \to T$ which sends $$(r,m) \times (r', m') \to (r, m+ m'),$$ where $t(r,m) = s(r',m')$ so that $(r,m) \times (r', m') \in T \times_{R, s, t} T.$
Therefore, we obtain a groupoid denoted as
$$M \times R  \, \substack{\longrightarrow \\ \longleftarrow \\ \longrightarrow}\, R. $$ Note that all the morphisms $s,t,c,e$ are morphisms of abelian groups, one can actually convert the above data into a $1$-truncated simplicial abelian group. 
\end{construction}{}
In the construction below, we explain how to attach a $1$-truncated simplicial ring or a ring groupoid from the data of a quasi-ideal.
\begin{construction}[Quasi-ideal as a simplicial commutative ring]\label{scr}Let $d: M \to R$ be a quasi-ideal. We have already defined a groupoid $$M \times R  \, \substack{\longrightarrow \\ \longleftarrow \\ \longrightarrow}\, R $$ which can also be thought of as a $1$-truncated simplicial abelian group. 
Next, we give a ring structure on $M \times R.$ 
We define $(m_1, r_1) \cdot (m_2, r_2) \coloneqq (r_2 m_1 + r_1 m_2 + d(m_1) m_2, r_1 r_2).$ Now as one easily checks, the morphisms $s,t,c,e$ in the definition of the groupoid $$M \times R  \, \substack{\longrightarrow \\ \longleftarrow \\ \longrightarrow}\, R $$ are all \textit{ring} homomorphisms with respect to the ring structure on $M \times R$ defined above. The above data can be converted into a $1$-truncated simplicial commutative ring.


\end{construction}{}

\end{definition}{}

\begin{definition}[Quasi-ideals in schemes]Let $R$ be a ring scheme and $M$ be a module scheme over $R$ equipped with a map $d: M \to R$ of $R$-module schemes. This data will be called a \textit{quasi-ideal in} $R$ if $d(x) \cdot y = d(y) \cdot x$ for scheme theoretic points $x,y \in M$.

A morphism between quasi-ideals in schemes is defined in a way similar to \cref{qid}.
\end{definition}{}

Finally, let us give some examples of quasi-ideals that will be used later on. For more details on these examples, we refer the reader to \cite[\S 3.2-3.5]{prismatization}
or \cite[\S 2.2]{Mon21}.

\begin{example}Let $\mathbb{G}_a^{\sharp} \to \mathbb{G}_a$ denote the quasi-ideal obtained by taking the divided power envelope of origin inside $\mathbb{G}_a$.
\end{example}{}

\begin{example}\label{gaperf}Let $B$ be any ring on which $p$ is nilpotent. Then the functor $S \to S^\flat:= \varprojlim_{F} S/p$ is representable by the affine ring scheme $\w{Spec}\, B[x^{1/p^\infty}]$ which will be denoted as $\mathbb{G}_a^{\w{perf}}.$

\end{example}{}

\begin{example}\label{gaperfsharp}Let $\mathbb{G}_a^{\w{perf} ,\sharp} \to \mathbb{G}_a^{\w{perf}}$ denote the quasi-ideal obtained by taking the divided power envelope of the closed subscheme defined by the ideal $(p,x)$ inside $\mathbb{G}_a^{\w{perf}}$ compatibly with the existing divided powers of $p.$
\end{example}{}

\begin{example} Let $W$ denote the ring scheme of $p$-typical Witt vectors. By taking kernel of the Frobenius $F,$ one obtains a quasi-ideal $W[F] \to \mathbb{G}_a,$ which is isomorphic to $\mathbb{G}_a^\sharp \to \mathbb{G}_a$ as a quasi-ideal in $\mathbb{G}_a.$
\end{example}{}

\begin{example}
By considering the multiplication by $p$ map on $W,$ one obtains a quasi-ideal $W \xrightarrow{\times p} W.$
\end{example}{}

\subsection{Ring stacks}
We begin by collecting some notations. If $C$ and $D$ denote two $\infty$-categories which have finite products, then the category of finite product preserving functors will be given by $\w{Fun}_{\times}({C}, {D}).$ Let $\w{Poly}_{A}$ denote the category of finitely generated polynomial algebras over $A.$

\begin{definition}[Animated ring objects in a category]\label{anim} Let $C$ be an $\infty$-category with products. Animated $A$-algebra objects in $C$, denoted as $\text{ARings}(C)_{A}$, is defined to be the category $\w{Fun}_{\times}(\w{Poly}_A^{\w{op}}, C)$.

In the case where $ C$ is the $\infty$-category of spaces, then the above definition with $A=\mathbb{Z}$ recovers the usual category of animated rings.
\end{definition}{}

\begin{remark}
The $\infty$-category of animated rings has all small colimits. Given a simplicial commutative ring, one can take colimit over the simplex category and obtain an animated ring. In particular, given a quasi-ideal, one can apply \cref{scr} and obtain an animated ring.
\end{remark}{}

\begin{definition}[Prestacks]\label{prestacks} The $\infty$-category of prestacks over a fixed (discrete) base ring $B$, denoted by $\w{PreSt}_B$, is defined to be the category of functors $\text{Fun} (\w{Alg}_B, \mathcal{S})$, where $\w{Alg}_B$ is the category of discrete $B$-algebras and $\mathcal{S}$ is the $\infty$-category of spaces.
\end{definition}{}

We note that even though we do not impose any sheafiness conditions, the examples of stacks we consider will all be (hypercomplete) fpqc sheaves of spaces.

\begin{definition}[$A$-algebra prestacks over $\mathrm{Spec}(B)$] 
The category of $A$-algebra prestacks over $\mathrm{Spec}(B)$, denoted as $A$-$\w{Alg-PreSt}_B$, 
is defined to be the category of animated $A$-algebra objects in the category $\w{PreSt}_B.$ 

\end{definition}{}

\begin{remark}
\label{switch}
We note that another possible way to define the category $A\w{-}\w{Alg}\w{-PreSt}_B$ is to define it as $\w{Fun}(\w{Alg}_B, \w{ARings}_{A}).$ However, this is equivalent to the definition considered above since we have natural equivalence of categories 
$$\w{Fun}(\w{Alg}_B, \w{ARings}_{A}) \simeq \w{Fun}(\w{Alg}_B, \w{Fun}_{\times}(\w{Poly}_A^{\w{op}}, \mathcal{S})) \simeq \w{Fun}_{\times}(\w{Poly}_A^{\w{op}}, \text{Fun} (\w{Alg}_B, \mathcal{S})) \simeq  \w{Fun}_{\times}(\w{Poly}_A^{\w{op}}, \w{PreSt}_B).$$ The middle equivalence uses the fact that product in functor category is calculated term-wise; the precise $\infty$-categorical (dual) assertion can be found in \cite[Corollary 5.1.2.3]{Lu}.
\end{remark}{}

\begin{construction}(Cone of a quasi-ideal) In view of \cref{switch} and \cref{scr}, it follows that given a quasi-ideal $d: M \to R$ in schemes, the quotient prestack $[R/M]$ (under the additive action of $M$ on the ring scheme $R$ by translation via $d$) has the structure of a ring prestack.
In the context of this paper, we will consider associated ring stacks of such ring prestacks, obtained by fpqc sheafification.

\end{construction}{}

\begin{example}
We will see later that all the examples of quasi-ideals from \cref{sec2.1} have the same cone.
\end{example}{}

\subsection{Affine stacks}

We will also use the notion of \textit{affine stacks} due to To\"{e}n \cite{MR2244263}. Here, we will recall its definition and basic properties very briefly, in the language of $\infty$-categories. To that end, we start by fixing an ordinary base ring $B$.
Let $\w{coSCR}_B$ denote the $\infty$-category of cosimplicial rings over $B$ arising from the simplicial model structure defined in \cite[Theorem 2.1.2]{MR2244263};  to construct the associated $\infty$-category from the simplicial model category, one looks at the fibrant simplicial category obtained from the subcategory of fibrant-cofibrant objects inside the given simplicial model category and applies the simplicial nerve construction, which produces an $\infty$-category by \cite[Proposition 1.1.5.10]{Lu}. It follows from \cite[Corollary 4.2.4.8]{Lu} that the $\infty$-category $\w{coSCR}_B$ has all small limits and colimits.
\begin{definition}[Affine stacks]\label{defaffinestack}
An object $\mathcal{Y}$ of $\w{PreSt}_B$ is called an \textit{affine stack} over $B$ if there exists an object $C \in \w{coSCR}_B$ such that $\mathcal{Y}$ is the restriction of the functor $h_C: \w{coSCR}_B \to \mathcal{S}$ corepresented by $C$ along the inclusion $\w{Alg}_B \to \w{coSCR}_B.$ The full subcategory of such objects inside $\w{PreSt}_B$ will be denoted by $\w{AffStacks}_B.$
\end{definition}{}

\begin{remark}\label{affcosimp} It follows from the definition that the category of affine stacks is stable under small limits, \textit{c.f.}~\cite[Proposition 2.2.7]{MR2244263}. Also, an affine stack is a hypercomplete fpqc sheaf of spaces, \textit{c.f.}~\cite[Lemma 1.1.2, Proposition 2.2.2]{MR2244263}. The key property of affine stacks that will be useful for us is the fact that taking derived global section functor induces an equivalence of $\infty$-categories $\w{AffStacks}_B \simeq \w{coSCR}_B^{\w{op}},$ \textit{c.f.}~\cite[Corollary 2.2.3]{MR2244263}. 
\end{remark}{}

\begin{remark}We point out that even though the definition of the subcategory of affine stacks $\w{AffStacks}_B$ inside $\w{PreSt}_B$ \textit{a priori} depends on the category $\w{coSCR}_B,$ the notion of being an affine stack is \textit{intrinsic}: being an affine stack is a property that can be formulated \textit{only} by using the fpqc topology and the category of ordinary rings. We refer to \cite[Theorem 2.2.9]{MR2244263} for a more precise formulation of this statement using Bousfield localization. \textit{A posteriori}, the same intrinsic property carries over to the $\infty$-category $\w{coSCR}_B$ which makes it rather special in comparison to certain other related categories such as the $\infty$-category of derived rings or $E_{\infty}$-rings.
\end{remark}{}

\begin{example} An affine scheme is clearly an affine stack. More precisely, the category $\w{Aff}_B$ of affine schemes over $B$ embeds fully faithfully inside the category $\w{AffStacks}_B$ of affine stacks over $B.$
\end{example}
\begin{example}
The stacks $K(\mathbb G_a, m)$ for $m \ge 0$ are examples of affine stacks \cite[Lemma 2.2.5]{MR2244263}. On the other hand, $K(\mathbb{G}_m, m)$ is \textit{not} an affine stack for any $m>0.$ By \cite[Corollary 2.4.10]{MR2244263}, for pointed and connected stacks over a field, being an affine stack is equivalent to the sheaf of all the higher homotopy groups being representable by unipotent affine group schemes (possibly of infinite type). 
\end{example}{}

\begin{remark}\label{hyper}
Let $\w{St}^\wedge_B$ denote the $\infty$-category of hypercomplete fpqc sheaves of spaces (see \cite[Section 6.5]{Lu} for a discussion of hypercomplete $\infty$-topos). Translating the results from \cite[Lemma 1.1.2, Proposition 2.2.2, Corollary 2.2.3]{MR2244263} in the language of $\infty$-categories, we obtain a colimit preserving functor $\w{St}^\wedge_B \to \w{coSCR}_B^{\w{op}}.$ There is also a natural colimit preserving functor $\w{PreSt}_B \to \w{St}^\wedge_B,$ and the composite functor denoted as $R\Gamma (\cdot, \mathcal{O})\colon \w{PreSt}_B \to \w{coSCR}_B^{\w{op}}$ gives us the ``derived global section functor.'' By construction, $R\Gamma(\cdot, \mathcal{O})\colon \w{PreSt}_B \to \w{coSCR}_B^{\w{op}}$ preserves all small colimits. By \cref{prestacks} and \cite[Lemma 5.1.5.5, Proposition 5.1.5.6]{Lu}, it follows that $R\Gamma(\cdot, \mathcal{O})$ can be simply described as left Kan extension of the functor $\w{Aff}_B \to \w{coSCR}_B^{\w{op}}$ (along the inclusion of categories $\w{Aff}_B \to \w{PreSt}_B$) which sends an affine scheme to its underlying ring of global sections. This checks the compatibility of two \textit{a priori} different ways of defining the derived global section functor.
\end{remark}
\begin{remark}\label{rightkan1}
Let us now suppose that $\mathcal{Y}$ is an affine stack over $B$ which is corepresented by $C \in \w{coSCR}_B$ (see \cref{defaffinestack}). As noted in \cref{affcosimp}, $\mathcal{Y}$ is a hypercomplete fpqc sheaf of spaces. According to \cite[Corollary 2.2.3]{MR2244263} and \cref{hyper}, we have a natural isomorphism $R\Gamma(\mathcal{Y}, \mathcal{O}) \simeq C$ in $\w{coSCR}_B.$ Unwrapping all the definitions and using the equivalence $\w{AffStacks}_B \simeq \w{coSCR}_B^{\w{op}},$ we obtain the categorical implication that the identity functor $\w{coSCR}_B \to \w{coSCR}_B$ is naturally equivalent to right Kan extension of the inclusion $\w{Alg}_B \to \w{coSCR}_B$ along itself. Roughly speaking, this means that for any $C \in \w{coSCR}_B,$ we have a natural isomorphism 
$$ C \simeq \left(\varprojlim_{\substack{C \to A,\\ \text{A is discrete}}} A \right) \in \w{coSCR}_B.$$
\end{remark}

\begin{remark}\label{rightkan}
The observation in \cref{rightkan1} regarding right Kan extension implies that if $\mathcal{D}$ is any $\infty$-category and $F: \w{coSCR}_B \to \mathcal{D}$ is a functor that is a right adjoint, then $F$ is naturally equivalent to the right Kan extension of the composite functor $\w{Alg}_B \to \w{coSCR}_B \to \mathcal{D}$ along the inclusion $\w{Alg}_B \to \w{coSCR}_B.$ 
\end{remark}{}

\subsection{Unwinding ring stacks}\label{unwindringstack}
In this section, we describe how to \textit{unwind} the data of a ring stack to obtain a cohomology theory. This construction is an $\infty$-categorical enhancement of \cite[Example 3.0.1]{Mon21} and we will call this the unwinding of a given ring stack. The construction only uses basic categorical principles such as Kan extensions and the necessary foundations can be found in \cite{Lu}.


\begin{construction}[Unwinding]
\label{unwinding construction}
We will construct a functor $$\mathrm{Un}: A\w{-Alg}\w{-PreSt}_B \to \w{Fun}(\w{ARings}_{A}, \w{CAlg}(D(B))).$$ Here, $\w{CAlg}(D(B))$ denotes the commutative algebra objects in the derived $\infty$-category $D(B).$ We think of the objects in the right hand side as ``algebraic cohomology theories''.

We begin by noting that by definition $A\w{-Alg}\w{-PreSt}_B \simeq \w{Fun}_{\times}(\w{Poly}_A^{\w{op}}, \w{PreSt}_B).$ By Kan extension, there is a derived global section functor $R\Gamma: \w{PreSt}_B \to \w{CAlg}(D(B))^{\w{op}}.$ By composition, we get a functor $$A\w{-Alg}\w{-PreSt}_B^{\w{op}} \to \w{Fun}(\w{Poly}_A, \w{CAlg}(D(B))).$$ Now, we can perform a left Kan extension along the inclusion $\w{Poly}_A \to \w{ARings}_{A}$ to obtain the desired unwinding functor 
$$ \mathrm{Un}: A\w{-Alg}\w{-PreSt}_B^{\w{op}} \to \w{Fun}(\w{ARings}_{A}, \w{CAlg}(D(B))).$$
\end{construction}{}

\begin{example}\label{gaidentity}When $A=B,$ and $\mathcal{Y} \in \w{PreSt}_B$ is taken to be the ring scheme $\mathbb{G}_{a,B},$ the functor $\mathrm{Un}(\mathbb{G}_{a,B})$ is simply the forgetful functor $\w{ARings}_{A} \to \w{CAlg}(D(B)).$
\end{example}{}

Below we will study compatibility of the unwinding construction with restriction of scalars. More precisely, let $\mathcal{Y} \in A\w{-Alg}\w{-PreSt}_B.$ Let $A' \to A$ be a map of discrete rings. Then there is an obvious functor $$\w{res}: A\w{-Alg}\w{-PreSt}_B \to   A'\w{-Alg}\w{-PreSt}_B.$$ Let $\mathcal{Y}':= \w{res}(\mathcal{Y}) \in A'\w{-Alg}\w{-PreSt}_B.$ Applying the unwinding construction, we obtain two functors $\w{Un}(\mathcal{Y}): \w{ARings}_{A} \to \w{CAlg}(D(B))$ and $\w{Un}(\mathcal{Y}'): \w{ARings}_{A'} \to \w{CAlg}(D(B)).$ Note that we also have a natural functor (given by the derived tensor product) $L: \w{ARings}_{A'} \to \w{ARings}_{A}.$ In this set up, we have the following compatibility.

\begin{proposition} 
\label{compatibility of unwinding}
We have $\mathrm{Un}(\mathcal{Y})\circ L \simeq \mathrm{Un}(\mathcal{Y}')$ in $\mathrm{Fun}(\mathrm{ARings}_{A}, \mathrm{CAlg}(D(B))).$
\end{proposition}{}

\begin{proof}Since $L$ is obtained by left Kan extension of the composite functor $\w{Poly}_{A'} \xrightarrow[]{\ell} \w{Poly}_A \xrightarrow[]{} \w{ARings}_{A},$ it would be enough to prove $\mathrm{Un}(\mathcal{Y})\circ \ell \simeq \mathrm{Un}(\mathcal{Y}')$ in $\mathrm{Fun}(\mathrm{Poly}_{A}, \mathrm{CAlg}(D(B)).$ By \cref{unwinding construction}, we note that $\mathcal{Y}$ is classified by a functor $U: \w{Poly}_A^{\w{op}} \to \w{PreSt}_B$ and $\mathcal{Y}'$ is classified by $U': \w{Poly}_{A'}^{\w{op}} \to \w{PreSt}_B,$ and for our purpose, it would be enough to prove that $U \circ \ell^{\w{op}} \simeq U'.$ By \cref{switch}, it would be enough to prove that the restriction of scalar functor $\w{ARings}_{A} \to \w{ARings}_{A'}$ is induced by $\ell^{\w{op}}$ under the identifications
$\w{ARings}_{A} \simeq \w{Fun}_{\times}(\w{Poly}_{A}^{\w{op}}, \mathcal{S})$ and $\w{ARings}_{A'} \simeq \w{Fun}_{\times}(\w{Poly}_{A'}^{\w{op}}, \mathcal{S}).$ But that follows from adjunction.\end{proof}{}

\begin{notation}\label{frobtwist} If $k$ is a perfect field of characteristic $p$ and $\mathcal{Y}$ is a $W_n(k)$-algebra stack for $1 \le n \le \infty$, 
then we will use $\mathcal{Y}^{(1)}$ to denote the $W_n(k)$-algebra stack obtained by restriction of scalars along the Witt vector Frobenius $W_n(k) \to W_n(k).$ \textit{c.f.}~\cref{compatibility of unwinding}.

\end{notation}{}

\begin{remark}
For this paper, the Frobenius twist $\mathcal{Y}^{(1)}$ of a stack $\mathcal{Y}$ will not play an important role because we always work over a perfect field and are interested in the question of endomorphisms of the stacks. Since it also does not change the underlying stack, for the most part of the paper, we will ignore this Frobenius twist.
\end{remark}{}

\begin{example}
The \Cref{compatibility of unwinding} shows that
the Frobenius twisted forgetful functor 
$$R \mapsto R^{(1)} \coloneqq R \otimes_{k, \mathrm{Frob}} k$$
from $\mathrm{ARings}_k \to \mathrm{CAlg}(D(k))$ is the unwinding of $\mathbb{G}_{a,k}^{(1)}.$ The relative Frobenius $R^{(1)} \to R$ can be obtained by undwinding the map of $k$-algebra stacks $\mathbb{G}_{a,k} \to \mathbb{G}_{a,k}^{(1)}$ induced by the Frobenius.
\end{example}{}

\subsection{de Rham cohomology via unwinding}
In this section, we will describe how to use the unwinding construction to recover de Rham or crystalline cohomology functors. To this end, let $n,m \ge 1$ be two arbitrary positive integers and let $p$ be a fixed prime. Further, we fix a perfect ring $k$ of characteristic $p.$ Let $W_r(k)$ denote the ring of $r$-truncated Witt vectors. Using crystalline cohomology, more precisely, its derived variant (see \cref{derivedcrys} below),
one obtains certain functors denoted as 
$$\mathrm{dR}_{m,n}: \w{ARings}_{W_n(k)} \to \w{CAlg}(D(W_m(k))),$$ 
which we will loosely still call de Rham cohomology functors and specify the $n,m$.
To define them one really needs to use a deformation of the de Rham cohomology functor, i.e.,~the crystalline cohomology functors. 

The following essentially already appeared in \cite[\S 10.2]{BLM20} and \cite[\S 8.2]{BMS2}.
\begin{definition}\label{derivedcrys}
Given a finitely generated polynomial $W_n(k)$-algebra $P$, define
$\mathrm{dR}_{m,n}(P) \coloneqq \mathrm{R\Gamma_{crys}}(P_0/W_m(k))$ where $P_0$ denotes the mod $p$ reduction of $P$.
We will let $\mathrm{dR}_{m,n}$ to denote the left Kan extension of the above functor from finitely generated polynomials to all animated $W_n(k)$-algebras which takes values in $\w{CAlg}(D(W_m(k))).$
\end{definition}
We use this notation as we believe that crystalline cohomology is secretly a disguise of derived
de Rham cohomology, see \cite[Proposition 3.27]{Bha12} and \cite[Proposition 2.11]{LL20} for occasions of this perspective.
Our goal is to describe $\mathrm{dR}_{m,n}$ as the unwinding of a certain object in $W_n(k) \w{-Alg-PreSt}_{W_m(k)}.$

\begin{definition}
\label{a1dR}
Let $W$ denote the ring scheme over $\mathrm{Spec}(W_m(k))$ underlying the $p$-typical Witt vectors. 
Using the Artin--Hasse homomorphism $W(k) \to W(W(k)),$ one can view $W$ as a $W(k)$-algebra scheme. Then $d: W^{(1)} \xrightarrow[]{\times p} W ^{(1)}$ defines a quasi-ideal in schemes. By considering its cone, one obtains a $k$-algebra stack over $\mathrm{Spec}(W_m(k))$, which can be regarded as a $W_n(k)$-algebra stack over $\mathrm{Spec}(W_m(k))$ via the natural map
$W_n(k) \twoheadrightarrow k$. We denote the resulting $W_n(k)$-algebra stack over $\mathrm{Spec}(W_m(k))$ by $\mathbb{A}^{1, \mathrm{dR}}_{m,n}.$ When $n$ is fixed, we will use $\mathbb{A}^{1, \w{dR}}_{B}$ to denote the pull-back of $\mathbb{A}^{1, \w{dR}}_{m,n}$ to $\w{Spec}\, B$ for a $W_m(k)$-algebra $B.$
\end{definition}{}

\begin{remark}The above definition gives a generalization of the definition of $\mathbb{A}^{1,\w{dR}}$ as an $\mathbb{F}_p$-algebra stack due to Drinfeld to the more general case of an arbitrary perfect ring $k.$ To do this, one crucially needs to  use the Artin--Hasse natural transformation $W(\cdot) \to W(W(\cdot))$. One can abstractly construct this natural transformation by realizing the functor $W$ as a right adjoint to the inclusion of the category of delta rings inside all rings.

\end{remark}{}

\begin{proposition}\label{affinestack}The stack underlying $\mathbb{A}^{1, \mathrm{dR}}_{m,n}$ is an affine stack.
\end{proposition}{}
 \begin{proof}Indeed, the stack underlying $\mathbb{A}^{1, \mathrm{dR}}_{m,n}$ is obtained by taking cone of $d: W \xrightarrow[]{\times p} W,$ which is the same as fibre of the induced map $BW \to BW.$ Since affine stacks are closed under limits, it would be enough to show that $BW$ is an affine stack. This follows from the proof of \cite[Proposition 3.2.7]{MRT19}.
 Let us give a rough sketch of their argument. Let $W_n$ denote the ring scheme underlying $n$-truncated $p$-typical Witt vectors. Using certain obstruction vanishing, one first argues that $B W \simeq \varprojlim B W_n.$ Therefore, it is enough to prove that $B W_n$ is an affine stack for all $n.$ 
 To do so, one argues by induction on $n$. Using the short exact sequence
 \[
 0 \to \mathbb{G}_a \to W_{n+1} \to W_n \to 0,
 \]
 one sees that $BW_{n+1}$ is classified by a map $BW_n \to K(\mathbb{G}_a, 2).$ More precisely, we have a fiber sequence 
 \begin{center}
     \begin{tikzcd}
BW_{n+1} \arrow[d] \arrow[rr] &  & * \arrow[d]         \\
BW_n \arrow[rr]               &  & {K(\mathbb{G}_a,2).}
\end{tikzcd}
 \end{center}{}
 Since the stacks $K(\mathbb{G}_a,m)$ are affine stacks for $m \ge 0,$ we are done by induction.
\end{proof}{}

\begin{remark} The above argument can be modified to more generally show that $K(W, m)$ is an affine stack for all $m\ge 0.$ Consequently, one can show that the abelian group stack $\mathbb{A}^{1, \w{dR}}[m]$ is also an affine stack for all $m \ge 0.$ We have $R\Gamma_{\w{dR}}(K(\mathbb{G}_a, m)) \simeq R\Gamma (\mathbb{A}^{1, \w{dR}}[m], \mathcal{O})$ for all $m \ge 0.$

\end{remark}{}

\begin{proposition}[{\textit{c.f.}~\cite[Remark 7.9]{BhaLur2}, ~\cite{Mon22}}]
\label{identify endomorphisms}
We have a natural isomorphism $\mathrm{Un}(\mathbb{A}^{1, \mathrm{dR}}_{m,n}) \simeq \mathrm{dR}_{m,n}.$
\end{proposition}{}


\begin{proof}By \cref{compatibility of unwinding}, the proof reduces to $n=1$. Further, by \cite[Theorem 1.1.1]{Mon21}, one can reduce to $m=1$. Let us now explain the proof of the natural isomorphism $\mathrm{Un}(\mathbb{A}^{1, \mathrm{dR}}_{1,1}) \simeq \mathrm{dR}_{1,1}.$
By construction of the unwinding functor, it would be enough to show that the restricted functors $\mathrm{dR}_{1,1} \colon \mathrm{Poly}_k \to \w{CAlg}(D(k))$ and $\mathrm{Un}(\mathbb{A}^{1, \mathrm{dR}}_{1,1}) \colon \mathrm{Poly}_k \to \w{CAlg}(D(k))$ are naturally isomorphic. Note that we have a natural functor $\w{coSCR}_k \to \w{CAlg}(D(k))$ of $\infty$-categories and by construction, $\mathrm{Un}(\mathbb{A}^{1, \mathrm{dR}}_{1,1})$ lifts to give a functor
still denoted as $\mathrm{Un}(\mathbb{A}^{1, \mathrm{dR}}_{1,1}): \mathrm{Poly}_k \to \w{coSCR}_k.$ By quasi-syntomic descent, $\mathrm{dR}_{1,1}$ also lifts to give a functor still denoted as $\mathrm{dR}_{1,1}: \mathrm{Poly}_k \to \w{coSCR}_k.$ It would be enough to prove that these two functors are naturally isomorphic.

By considering $\w{gr}^0$ of the Hodge filtration on de Rham cohomology and quasi-syntomic descent, there is a natural arrow $\w{dR}_{1,1} \to \iota$ in the category $\w{Fun}(\w{Poly}_k, \w{coSCR}_k),$ where $\iota: \w{Poly}_k \to  \w{coSCR}_k$ denotes the natural inclusion functor.

Note that the derived global sections of $\mathbb{A}^{1,\w{dR}}_{1,1}$ agrees with $\w{dR}_{1,1}(k[x])$ in $\w{coSCR}_k.$ For this, one can use the identification $\w{Cone}(\mathbb{G}_a^\sharp \to \mathbb{G}_a) \simeq \w{Cone}(W \xrightarrow[]{\times p} W)$ and the \v{C}ech--Alexander complex. Since by \cref{affinestack}, $\mathbb{A}^{1, \w{dR}}_{1,1}$ is an affine stack, it follows that the functors $\mathrm{Un}(\mathbb{A}^{1, \mathrm{dR}}_{1,1}): \mathrm{Poly}_k \to \w{coSCR}_k$ and $\mathrm{dR}_{1,1}: \mathrm{Poly}_k \to \w{coSCR}_k$ preserve finite coproducts. In order to check that they are naturally isomorphic, it is enough to do so for the functors $\mathrm{Un}(\mathbb{A}^{1, \mathrm{dR}}_{1,1})': \mathrm{ARings}_k \to \w{coSCR}_k$ and $\mathrm{dR}_{1,1}': \mathrm{ARings}_k \to \w{coSCR}_k$ obtained by left Kan extension along $\w{Poly}_k \to \w{ARings}_k.$ 

By \cite[Proposition 5.5.8.15]{Lu}, the functors $\mathrm{Un}(\mathbb{A}^{1, \mathrm{dR}}_{1,1})'$ and $\mathrm{dR}_{1,1}'$ both preserve small colimits. Similarly, by left Kan extension, $\iota: \w{Poly}_k \to  \w{coSCR}_k$ also extends to a colimit preserving functor $\iota': \w{ARings}_k \to \w{coSCR}_k.$ By the adjoint functor theorem, all of these functors have right adjoints. Let ${\mathrm{Un}(\mathbb{A}^{1, \mathrm{dR}}_{1,1})'}{^\mathcal{R}}, {\mathrm{dR}_{1,1}'}{^\mathcal{R}}$ and ${\iota'} ^{\mathcal{R}}$ denote the right adjoints to $\mathrm{Un}(\mathbb{A}^{1, \mathrm{dR}}_{1,1})', \mathrm{dR}_{1,1}'$ and $\iota'$ respectively. It would be enough to prove that ${\mathrm{Un}(\mathbb{A}^{1, \mathrm{dR}}_{1,1})'}{^\mathcal{R}} \simeq {\mathrm{dR}_{1,1}'}{^\mathcal{R}}.$

Let ${\mathrm{Un}(\mathbb{A}^{1, \mathrm{dR}}_{1,1})'}{^\mathcal{R}_\circ}, {\mathrm{dR}_{1,1}'}{^\mathcal{R}_\circ}$ and ${\iota'} {^\mathcal{R}_\circ}$ denote the restriction of the functors ${\mathrm{Un}(\mathbb{A}^{1, \mathrm{dR}}_{1,1})'}{^\mathcal{R}}, {\mathrm{dR}_{1,1}'}{^\mathcal{R}}$ and ${\iota'} ^{\mathcal{R}}$ along the inclusion of categories $\w{Alg}_k \to \w{coSCR}_k$ respectively. For our purpose, by considering right Kan extensions as explained in \cref{rightkan}, it would be enough to prove that ${\mathrm{Un}(\mathbb{A}^{1, \mathrm{dR}}_{1,1})'}{^\mathcal{R}_\circ} \simeq {\mathrm{dR}_{1,1}'}{^\mathcal{R}_\circ}.$ Note that they are both functors from $\w{Alg}_k$ to $\w{ARings}_k.$ Further, one observes that for an $S \in \w{Alg}_k,$ we have ${\iota'} {^\mathcal{R}_\circ}(S)= S,$ which identifies with the $S$-valued points of the ring scheme $\mathbb{G}_a.$ Thus we have a natural arrow ${\iota'} {^\mathcal{R}_\circ} \simeq \mathbb{G}_a \to {\mathrm{dR}_{1,1}'}{^\mathcal{R}_\circ}$ in $\w{Fun}(\w{Alg}_k, \w{ARings}_k),$ where $\mathbb{G}_a$ is viewed as an object of $\w{Fun}(\w{Alg}_k, \w{ARings}_k)$ by considering its functor of points. We note the following lemma.

\begin{lemma}
The fiber $F$ of the map $\mathbb{G}_a \to {\mathrm{dR}_{1,1}'}{^\mathcal{R}_\circ}$ identifies with the $\mathbb{G}_a$-module scheme $\mathbb{G}_a^\sharp.$ 
\end{lemma}
\begin{proof}To see this, we note that ${\mathrm{dR}_{1,1}'}{^\mathcal{R}_\circ}$ can be viewed as a ring stack, whose underlying stack, by construction, is the affine stack corresponding to the object $\w{dR}'_{1,1} (k[x]) \in \w{coSCR}_k.$ Therefore, the stack underlying $F$ is given by the affine stack corresponding to the cosimplicial ring obtained by the pushout $\displaystyle{k \sqcup_{\w{dR}'_{1,1}(k[x])} k[x]}$ in $\w{coSCR}_k.$ Since $\mathbb{A}^{1, \w{dR}}_{1,1}$ is an affine stack and $R\Gamma (\mathbb{A}^{1, \w{dR}}_{1,1}, \mathcal{O}) \simeq \w{dR}'_{1,1}(k[x]),$ it follows that $k \sqcup_{\w{dR}'_{1,1}(k[x])} k[x] \simeq R\Gamma(\mathbb{G}_a^\sharp, \mathcal{O}) \simeq D_{x}(k[x]),$ where the latter denotes the divided power envelope of $k[x]$ at the ideal $(x).$ In particular, the pushout is a discrete ring and the stack underlying $F$ is an affine scheme. Let $\w{dR}(k[x])$ denote the object of $\w{CAlg}(D(k))$ underlying $\w{dR}'_{1,1}(k[x]).$ Then there is a natural map $$k \otimes_{\w{dR}(k[x])}k[x] \to k \sqcup_{\w{dR}'_{1,1}(k[x])} k[x]$$ in $\w{CAlg}(D(k)).$ We have an isomorphism $k \otimes_{\w{dR}(k[x])}k[x] \simeq \w{dR}_{k/k[x]},$ where the latter denotes derived de Rham cohomology. By \cite[Lemma 3.29]{Bha12}, it follows that $\w{dR}_{k/k[x]} \simeq D_x(k[x])$
and the natural map above is an isomorphism. We note that $\w{Spec}(k \otimes_{\w{dR}(k[x])}k[x]) \simeq \w{Spec}(\w{dR}_{k/k[x]})$ also has the structure of a group scheme, where the multiplication is induced by functoriality of $\w{dR}_{k/ (\cdot)}$
along the map $k[x] \to k[x] \otimes_k k[x]$ given by $x \to x \otimes 1 + 1 \otimes x.$ Moreover, $\w{Spec}(\w{dR}_{k/k[x]})$ has the structure of a $\mathbb{G}_a$-equivariant group scheme, where the $\mathbb{G}_a$-action is given by the map $k[x] \simeq \w{dR}_{k[x]/k[x]} \to \w{dR}_{k/k[x]} \otimes_k k[x] \simeq D_{x}(k[x]) \otimes_k k[x]$ which is induced by functoriality of derived de Rham cohomology applied to the diagram below.

\begin{center}
\begin{tikzcd}
{k[x]} \arrow[rr, "x \to x"] \arrow[d, "x \to x\otimes x"']           &  & k[x] \arrow[d, "x \to 0"] \\
{k[x]\otimes k[x]} \arrow[rr, "{x\otimes 1 \to 0, 1\otimes x\to x}"'] &  & {k[x]}     
\end{tikzcd}
\end{center}{}
Using the explicit description of the induced maps, one explicitly verifies that $\w{Spec}(\w{dR}_{k/k[x]})$ is naturally isomorphic to $\mathbb{G}_a^\sharp$ as a $\mathbb{G}_a$-module scheme. Further, by applying functoriality along the diagrams mentioned earlier, we see that the map of schemes $F \to \w{Spec}(\w{dR}_{k/k[x]})$ induced by the natural map $k \otimes_{\w{dR}(k[x])}k[x] \to k \sqcup_{\w{dR}'_{1,1}(k[x])} k[x]$ above is actually a $\mathbb{G}_a$-equivariant map of group schemes. Since we have already noted that $k \otimes_{\w{dR}(k[x])}k[x] \to k \sqcup_{\w{dR}'_{1,1}(k[x])} k[x]$ is an isomorphism, this shows that $F$ is indeed isomorphic to $\mathbb{G}_a^\sharp$ as a $\mathbb{G}_a$-module scheme, as desired.
\end{proof}

Now we have obtained a natural map $\mathbb{A}^{1, \w{dR}}_{1,1} \simeq \w{Cone}(\mathbb{G}_a^\sharp \to \mathbb{G}_a) \to {\mathrm{dR}_{1,1}'}{^\mathcal{R}_\circ}$ of $k$-algebra stacks, i.e., as objects in the category $\w{Fun}(\w{Alg}_k, \w{ARings}_k).$ We have already noted that their underlying stacks are isomorphic. Thus we obtain an isomorphism  
$\mathbb{A}^{1, \w{dR}}_{1,1} \simeq {\mathrm{dR}_{1,1}'}{^\mathcal{R}_\circ}.$
Since the stack underlying $\mathbb{A}^{1, \w{dR}}_{1,1}$ is an affine stack, it follows that ${\mathrm{Un}(\mathbb{A}^{1, \mathrm{dR}}_{1,1})'}{^\mathcal{R}_\circ} \simeq \mathbb{A}^{1, \w{dR}}_{1,1}$ as objects of $\w{Fun}(\w{Alg}_k, \w{ARings}_k).$ This constructs the isomorphism ${\mathrm{Un}(\mathbb{A}^{1, \mathrm{dR}}_{1,1})'}{^\mathcal{R}_\circ} \simeq  {\mathrm{dR}_{1,1}'}{^\mathcal{R}_\circ},$ which finishes the proof.
\end{proof}{}
The following fact was used in the above proof, which uses compatibility of two models of the $k$-algebra stack $\mathbb{A}^{1, \mathrm{dR}}$ over $\mathrm{Spec}(W(k))$.
\begin{proposition}[{\cite[3.5.1]{prismatization}}]
\label{two models}
There is an isomorphism of $k$-algebra stacks over $\mathrm{Spec}(W(k))$:
\[
\mathrm{Cone}(\mathbb{G}_a^{\sharp} \to \mathbb{G}_a) \simeq \mathrm{Cone}(W^{(1)} \xrightarrow{\times p} W^{(1)}).
\]
\end{proposition}

The $k$-algebra structure on the source comes from
the natural maps $W(k) \to \mathbb{G}_a$ and $W(k) \xrightarrow{1 \mapsto p - V(1)} W[F]$.
To see the two underlying abelian group stacks are the same, 
notice that we always have $FV = p$ on the $p$-typical Witt ring, hence we get a factorization
\[
\xymatrix{
W^{(1)} \ar[rr]^{\times p} \ar[rd]^{V} & & W^{(1)}. \\
& W \ar[ru]^{F} &
}
\]
One then applies the octahedral axiom to the above triangle.
The fact that it induces an algebra isomorphism
can be seen using the fact that $F$ is an algebra homomorphism.
Said differently, one pulls back the quasi-ideal $W^{(1)} \xrightarrow{\times p} W^{(1)}$
along $W \xrightarrow{F} W^{(1)}$ to build the intermediate model relating the above two models.

\vspace{2 mm}

\begin{remark}\label{hodgeconj}
Note that there is a natural map of $k$-algebra stacks $\mathbb{G}_a \to \mathbb{A}^{1, \mathrm{dR}}$ whose unwinding provides a natural transformation $\mathrm{dR}(S) \to S,$ which corresponds to the natural projection onto the $\mathrm{gr}^{0}$ of the Hodge filtration on de Rham cohomology. There is also a natural map $\mathbb{A}^{1, \mathrm{dR}} \to \pi_0 (\mathbb{A}^{1, \mathrm{dR}}) = \mathbb{G}_a^{(1)}$ of $k$-algebra stacks which unwinds to the natural transformation $S^{(1)} \to \mathrm{dR}(S)$ induced by $\mathrm{Fil}^{0}$ of the conjugate filtration, \textit{c.f.}~\cref{cohomology of A1dR}.
\end{remark}{}

Now, we will see that the quasi-ideal $\mathbb{G}_a^{\w{perf}, \sharp} \to \mathbb{G}_a^{\w{perf}}$ that appears in \cite[Proposition 4.0.11]{Mon21} gives a third model of the $k$-algebra stack $\mathbb{A}^{1, \mathrm{dR}}$ over $\mathrm{Spec}(W(k))$, see also \cite{drinfeld2018stacky}. First, we will make some preparations. Below we always fix a positive integer $m$.

\begin{lemma}
\label{Wperf lemma}
On the fpqc site of $\mathbb{Z}/p^m$, we have 
$R\varprojlim_F W \simeq \varprojlim_F W$, which is representable by an affine scheme.
Moreover, its functor of points can be described as
$B \mapsto W(B^{\flat})$.

We denote the affine scheme representing $\varprojlim_F W$ by $W^{\mathrm{perf}}$. This scheme can be given an $W(k)$-algebra scheme structure when viewed over $\mathrm{Spec}(W_m(k))$.
\end{lemma}
\begin{proof}
The first assertion follows from \cite[Example 3.1.7 and Proposition 3.1.10]{proetale} and the fact that $F$ on $W$
is faithfully flat.
Inverse limit of affine schemes is again affine. For the last claim, we consider the following diagram of fpqc sheaves as a pro-object.
\[
\xymatrix{
\ldots \ar[r] & W_3 \ar[r]^{F} & W_2 \ar[r]^{F} & W_1 \\
\ldots \ar[r] & W_4 \ar[r]^{F} \ar[u]^{R} & W_3 \ar[r]^{F} \ar[u]^{R} & W_2 \ar[u]^{R} \\
\ldots \ar[r] & W_5 \ar[r]^{F} \ar[u]^{R} & W_4 \ar[r]^{F} \ar[u]^{R} & W_3 \ar[u]^{R} \\
\vdots \ar[r] & \vdots \ar[r]^{F} \ar[u]^{R} & \vdots \ar[r]^{F} \ar[u]^{R} & \vdots \ar[u]^{R}
}
\]
Taking limit vertically and then horizontally gives us $\varprojlim_F W$. Next we take limit horizontally and then vertically instead. Taking limits horizontally, we obtain the sheaf that sends $B \mapsto \varprojlim_F W_r(B)$,
which is canonically identified with $W(B^{\flat})$ by \cite[Lemma 3.2]{BMS1} (with $\pi$ in loc.~cit.~being $p$).
The vertical map $R$ is actually an isomorphism now, by the same \cite[Lemma 3.2]{BMS1}. This gives $\varprojlim_F W (B) \simeq W(B^\flat),$ as desired.
\end{proof}

Recall that the $F$ on $W$ induces a map $\mathbb{A}^{1, \mathrm{dR}} \to \mathrm{Frob}_{k,*}\mathbb{A}^{1, \mathrm{dR}}$
of $k$-algebra stacks which we will again denote by $F$.
We may untwist the Frobenius using inverse of Frobenius on $k$ on the source of this map.
Therefore, we get a $k$-algebra structure on the stack $\varprojlim_F(\mathbb{A}^{1, \mathrm{dR}})$.

\begin{lemma}\label{perfectionlemma}
We have an isomorphism of $k$-algebra stacks
$\mathbb{G}_a^{\mathrm{perf}} \simeq \varprojlim_F(\mathbb{A}^{1, \mathrm{dR}})$ over $\mathrm{Spec}(W_m(k))$.
\end{lemma}

\begin{proof}
By \Cref{Wperf lemma}, we see that 
\[
        R\varprojlim_F(\mathbb{A}^{1, \mathrm{dR}}) = \mathrm{Cone}({R}\varprojlim_F W \xrightarrow{\times p} {R}\varprojlim_F W)
= \mathrm{Cone}(W^{\mathrm{perf}} \xrightarrow{\times p} W^{\mathrm{perf}})
\]
whose functor of points is given by $B \mapsto B^{\flat}$.
Hence $\varprojlim_F(\mathbb{A}^{1, \mathrm{dR}})$ is isomorphic to $\mathbb{G}_a^{\mathrm{perf}}$ as a $k$-algebra stack (and in fact is a scheme).
\end{proof}

Therefore, we get a map of $k$-algebra stacks $\mathbb{G}_a^{\mathrm{perf}} \to \mathbb{A}^{1, \mathrm{dR}}$ 
over $\mathrm{Spec}(W_m(k))$.

\begin{lemma}
\label{perf model cover lemma}
The map of ($k$-algebra) stacks $f \colon \mathbb{G}_a^{\mathrm{perf}} \to \mathbb{A}^{1, \mathrm{dR}}$ is faithfully flat.
\end{lemma}

\begin{proof}
We look at the following diagram of $k$-algebra stacks
\[
\xymatrix{
W^{\mathrm{perf}} \ar[rr] \ar[rd] & & \mathbb{G}_a^{\mathrm{perf}} \ar[ld] \\
& \mathbb{A}^{1, \mathrm{dR}} &
}
\]
and observe that the horizontal and the left arrow are faithfully flat, hence the right arrow is faithfully flat as well.
\end{proof}

Let $K$ be the quasi-ideal in $\mathbb{G}_a^{\mathrm{perf}}$ given by the kernel of $f$,
then \Cref{perf model cover lemma} implies that $f$ gives rise to an isomorphism
of $k$-algebra stacks
$\mathrm{Cone}(K \to \mathbb{G}_a^{\mathrm{perf}}) \simeq \mathbb{A}^{1, \mathrm{dR}}$.
This is would be what we called the third model of $\mathbb{A}^{1, \mathrm{dR}}$; to complete its description, it remains
to understand the quasi-ideal $K$.

\begin{proposition}
$K$ is isomorphic to $\mathbb{G}_a^{\mathrm{perf}, \sharp}$ as quasi-ideals in $\mathbb{G}_a^{\mathrm{perf}}$.
In particular, $\mathrm{Cone}(\mathbb{G}_a^{\mathrm{perf}, \sharp} \to \mathbb{G}_a^{\mathrm{perf}}) \simeq \mathbb{A}^{1, \mathrm{dR}}$ as $k$--algebra stacks.
\end{proposition}

\begin{proof}This assertion follows from applying the (derived) crystalline cohomology functor $R\Gamma_{crys}$ to the following pushout diagram 
\begin{center}
\begin{tikzcd}
k \arrow[r]                               & {k[x^{1/p^\infty}]/x}         \\
{k[x]} \arrow[r] \arrow[u, "x \mapsto 0"] & {k[x^{1/p^\infty}]} \arrow[u]
\end{tikzcd}
\end{center}
and noting that global sections of $\mathbb{G}_a^{\w{perf}, \sharp}$ recovers $R\Gamma_{\w{crys}} (k[x^{1/p^\infty}]/x)$ and $R\Gamma_{\w{crys}}$ preserves the pushout diagram.
\end{proof}{}

\begin{remark} Using the above methods, let us sketch a quick proof of a result due to Bhatt--Lurie--Mathew \cite[Proposition 10.3.1]{BLM20}, \textit{c.f.}~\cite[Proposition 4.0.7]{Mon21}. Under \cref{identify endomorphisms}, \cref{gaidentity} and \cref{hodgeconj}, the assertion amounts to studying endomorphisms of $\mathbb{A}^{1,\w{dR}}$ respecting the natural map $\mathbb{G}_a \to \mathbb{A}^{1,\w{dR}}.$ By \cref{two models}, it is enough to show that the quasi-ideal $\mathbb{G}_a^{\sharp} \to \mathbb{G}_a$ has no non-trivial endomorphism as a quasi-ideal in $\mathbb{G}_a.$ But this follows directly from graded Cartier duality \cite[\S 2.4]{Mon21}.

\end{remark}{}

\begin{remark} Note that the definition of $\mathbb{A}^{1,\mathrm{dR}}$ as a $k$-algebra stack differs from $\mathrm{Cone}(W \xrightarrow[]{\times p} W)$ by a Frobenius twist. Indeed, the latter unwinds to Hodge-Tate cohomology (or a suitable base change of prismatic cohomology) \cite{BS19} which is the Frobenius descent
of de Rham cohomology (or crystalline cohomology)
\textit{c.f.}~\cref{compatibility of unwinding}.
\end{remark}{}

\section{Endomorphisms of de Rham cohomology I}
\label{section 3}

The quasi-syntomic descent technique introduced in \cite{BMS2} is a powerful tool in calculating
endomorphisms of de Rham cohomology functors in various settings.
We will illustrate them in this section.

Let $A \to B$ be a map of derived $p$-complete rings with bounded $p^\infty$-torsion.
In this section we consider the functor that sends any derived $p$-complete $A$-algebra $R$ to
\[
(\mathrm{dR} \widehat{\otimes}_{A} B)(R) \coloneqq \mathrm{dR}_{R/A} \widehat{\otimes}_{A} B \in \w{CAlg}(D(B)),
\]
where $\mathrm{dR}_{R/A}$ denotes the $p$-adic derived de Rham complex of $R$ relative to $A$
and $\widehat{\otimes}$ denotes derived $p$-completed tensor product.
If $B = A$, we simply denote the functor by $\mathrm{dR}$.

We are interested in the space of endomorphisms of this functor, viewed (by left Kan extension) as an object in the $\infty$-category
of functors from the $\infty$-category of derived $p$-complete animated rings to $\w{CAlg}(D(B))$. Let $\mathrm{qSyn}_A$ denote small quasi-syntomic site of $A$ which consists of algebras that are quasi-syntomic over $A$ and the covers are given by quasi-syntomic covers (see \cite[\S 4.2]{BMS2}).

\begin{proposition}[{\textit{c.f.}~\cite[Example 5.12]{BMS2}}]
\label{p-completed base changed dR is a qsyn sheaf}
The functor $\mathrm{dR} \widehat{\otimes}_{A} B$ when restricted to $\mathrm{qSyn}_A$ defines a quasi-syntomic sheaf.
\end{proposition}

\begin{proof}
It suffices to check this after going derived modulo $p$, so we are reduced to checking the following:
given $R \to S$ a faithfully flat quasi-syntomic map of algebras in $\mathrm{qSyn}_A$
with its \v{C}ech nerve $S^{\bullet}$,
then there is an isomorphism:
\[
\mathrm{dR}_{R/A} \otimes_{A} B/p \simeq \lim\left(\mathrm{dR}_{S^{\bullet}/A} \otimes_{A} B/p\right),
\]
where $B/p$ is the animated ring $B \otimes_{\mathbb{Z}} \mathbb{F}_p$.
By base change of derived de Rham cohomology, proving the above is equivalent to showing that:
\[
\mathrm{dR}_{(R \otimes_{A} B/p)/(B/p)} \simeq \lim\left(\mathrm{dR}_{(S^{\bullet} \otimes_{A} B/p)/(B/p)}\right),
\]
see \cite[p.33-35]{KP21} for a discussion of derived de Rham cohomology of maps of
animated rings.
To prove the above isomorphism, we employ the conjugate filtration \cite[Construction 2.3.12]{KP21}
(with base ring being $\mathbb{F}_p$).
The conjugate filtration is exhaustive and uniformly bounded above $-1$, hence it suffices to prove
its graded pieces satisfy similar quasi-syntomic descent.
Using the description of graded pieces of conjugate filtration we are finally reduced to showing that:
\[
\wedge^i_A \mathbb{L}_{R/A} \otimes_A \varphi_*(B/p) \simeq 
\lim \left(\wedge^i_{S^{\bullet}} \mathbb{L}_{S^{\bullet}/A} \otimes_{A} \varphi_*(B/p)\right).
\]
Here $\varphi_*(B/p)$ expresses the $A$-module structure on $B/p$ which is given by
$A \to B \to B/p \xrightarrow{\varphi} B/p$.
The next \Cref{descent of base changed cotangent complex} finishes the proof.
\end{proof}

\begin{proposition}[Flat descent for ``tensored'' cotangent complex]
\label{descent of base changed cotangent complex}
Fix a base ring $A$. For each $n \geq 0$ and an object $M \in D(A)$,
the functor $R \mapsto \wedge^n_R \mathbb{L}_{R/A} \otimes_A M$ is an fpqc sheaf with values in the $\infty$-category $D(A).$
\end{proposition}

\begin{proof}
One simply runs through the proof of \cite[Theorem 3.1]{BMS2} and sees that it works in this generality.
For convenience of the reader, let us illustrate the proof when $n = 1$.
Let $R \to S$ be a faithfully flat map of $A$-algebras with \v{C}ech nerve $S^{\bullet}$.
Using the transitivity triangle associated to 
$A \to R \to S^{\bullet}$ and applying the exact functor $(-) \otimes_A M$, we get a cosimplicial exact triangle
\[
\mathbb{L}_{R/A} \otimes_R S^{\bullet} \otimes_A M \to \mathbb{L}_{S^{\bullet}/A} \otimes_A M \to \mathbb{L}_{S^{\bullet}/R} \otimes_A M.
\]
We are therefore reduced to showing that:
\begin{itemize}
\item The map $R \to S^{\bullet}$ induces an isomorphism $\mathbb{L}_{R/A} \otimes_A M \to \lim \mathbb{L}_{R/A} \otimes_A M \otimes_R S^{\bullet}$;
\item $\lim \mathbb{L}_{S^{\bullet}/R} \otimes_A M = 0$.
\end{itemize}
The first item follows from fpqc descent along $R \to S$ by considering $\mathbb{L}_{R/A} \otimes_A M \in D(R).$ The second item is proved via a few reduction steps. By the convergence of the Postnikov filtration, it is enough to show that $\lim \pi_i(\mathbb{L}_{S^\bullet/R} \otimes_A M) \simeq 0$ in $D(R)$ for an arbitrary $i \in \mathbb{Z}$ which will be fixed from now. Again, by faithfully flat descent, it suffices to check $\left (\lim \pi_i(\mathbb{L}_{S^\bullet/R} \otimes_A M)\right )\otimes_R S \simeq \lim \left( \pi_i(\mathbb{L}_{S^\bullet/R} \otimes_A M)\otimes_R S \right) \simeq  \lim \pi_i(\mathbb{L}_{S^\bullet/R}\otimes_R S \otimes_A M ) \simeq 0.$ Let $S \to T^\bullet$ denote the base change of $R \to S^\bullet$ along $R \to S.$ By base change for cotangent complex, we need to show that $\lim \pi_i(\mathbb{L}_{T^\bullet/S} \otimes_A M ) \simeq 0.$ Since $S \to T^\bullet$ is the \v{C}ech nerve of the map $S \to S \otimes_R S,$ which admits a section, it follows that $S \to T^\bullet$ is a homotopy equivalence of cosimplicial $S$-algebras. Now we observe that $F:=\pi_i (\mathbb{L}_{(\cdot)/S} \otimes_A M)$ is a functor from the category of $S$-algebras to the category of abelian groups. Therefore, the cosimplicial abelian group $F(T^\bullet)$ is homotopy equivalent to $F(S).$ Since $F(S) \simeq 0,$ we obtain $\lim \pi_i(\mathbb{L}_{T^\bullet/S} \otimes_A M ) \simeq 0,$ as desired.
\end{proof}

As a consequence, let us record a result that says that the space of endomorphisms is actually discrete,
i.e., the homotopy groups in degrees above zero are trivial for every choice of base points.

\begin{lemma}
\label{discrete endo space}
The space of endomorphisms $\mathrm{End}(\mathrm{dR} \widehat{\otimes}_{A} B)$ is discrete.
\end{lemma}

\begin{proof}
First observe that $\mathrm{dR} \widehat{\otimes}_{A} B$ is left Kan extended from its restriction to the category of $p$-completely
finitely generated polynomial $A$-algebras.
Hence the restricted functor has the same space of endomorphisms.
Since our functor $\mathrm{dR} \widehat{\otimes}_{A} B$ is a sheaf 
on the quasi-syntomic site of $A$
and since $p$-completed polynomial $A$-algebras are quasi-syntomic over $A$,
restricting our functor to the full subcategory of $A$-algebras consisting of algebras that are quasi-syntomic over $A$
again computes the same endomorphism space.
Recall that since the quasi-syntomic site of $A$ admits a basis consisting of large quasi-syntomic $A$-algebras (see \cite[Definition 15.1]{BS19}), we may restrict our (base-changed) de Rham cohomology functor
to this basis and compute the space of endomorphisms there.
But now the values of the de Rham cohomology functor are $p$-completely flat $A$-algebras,
hence the base-changed de Rham cohomology functor has values which are discrete $B$-algebras
\cite[Lemma 4.6]{BMS2},
consequently the space of endomorphisms is discrete.
\end{proof}

If $\w{Spf}(A)$ has a disconnection, then the space of endomorphisms will be the product of endomorphism spaces on each
subset giving rise to the disconnection.
Hence without loss of generality, let us only treat those $A$'s with connected formal spectrum.
The following simple lemma will be used later, so let us record it here.

\begin{lemma}
\label{F_q lemma}
Let $A_0$ be an idempotent-free $\mathbb{F}_p$-algebra. Let $q$ be a power of $p$.
If every element $a \in A_0$ satisfies $a^q = a$, then $A_0$ is a sub-field inside $\mathbb{F}_q$.
\end{lemma}

In the rest of this section we will compute the space of endomorphisms in two cases

\vspace{1mm}
\noindent
Case I: When $A$ is the Witt ring of an idempotent-free characteristic $p$ perfect ring $k$ and $B = A$; and \vspace{1mm}

\noindent
Case II: When $A$ is a perfect $\mathbb{F}_p$-algebra and $B$ is an arbitrary $A$-algebra.
\vspace{1mm}

Building on the method of \cite[10.3-10.4]{BLM20},
Case I above is essentially worked out in the proof of \cite[Theorem 3.13]{LL20}, let us state a slightly more general result below.
\begin{proposition}
\label{integral endomorphism prop}
Assume that $A$ is $p$-torsion free, $p$-adically complete and $\mathrm{Spec}(A/p)$ is reduced and connected.
Then
$$
\mathrm{End}(\mathrm{dR}) =
\begin{cases}
\mathrm{Frob}_q^{\mathbb{N}} & \text{if}\, A = \mathbb{Z}_q \coloneqq W(\mathbb{F}_q), \\
\mathrm{id} & \text{otherwise}.
\end{cases}
$$
\end{proposition}

In \cite[\S 2.3]{LL20}, a Frobenius map is constructed on $p$-adic derived de Rham cohomology when the base
is a $p$-torsion free $\delta$-ring and it is semi-linear with respect to the Frobenius on the base $\delta$-ring.
The $\mathrm{Frob}_q$ appearing above is the corresponding power of the Frobenius
associated with the base $\delta$-ring $\mathbb{Z}_q$, one checks
easily that it is $\mathbb{Z}_q$-linear as desired.

\begin{proof}
Let us use $\w{Perf}$ to denote the full subcategory of those $A$-algebras
which are of the form $A\langle X_h^{1/p^\infty} \mid h \in H\rangle$ where $H$ is a set.
The proof of \cite[Theorem 3.13]{LL20} shows that
\begin{itemize}
    \item By restricting our de Rham cohomology functor to $\w{Perf}$, we get an injection of endomorphism monoids;
    \item The restricted de Rham cohomology functor has endomorphism monoid given by a sub-monoid in $\mathbb{Z}$;
    \item An element $n \in \mathbb{Z}$ above is characterized by its effect on $R = A\langle X^{1/p^\infty} \rangle$,
    which sends $X \mapsto X^{p^n}$; and
    \item The image of the restriction map is contained in $\mathbb{N} \subset \mathbb{Z}$.
\end{itemize}

Let us assume that $q = p^n$ is in the image of the restriction map.
Let $R = A\langle X^{1/p^\infty}\rangle$.
Take any $a \in A$ and let us contemplate the map $R \to R/(X-a)$.
The induced map of de Rham cohomology is the natural inclusion $R = A\langle X^{1/p^\infty} \rangle \to D$,
where $D$ is the algebra obtained by $p$-completely adjoining divided powers of $X-a$ to $R$.
To extend the map $X \mapsto X^{q}$ from $R$ to $D$, it is the same as requiring the image of $X-a$
to have divided powers.
Since in $D/p$ we have $X^p = a^p$, we see that $X^q - a = a^q - a + p \cdot d$ for some $d \in D$.
The condition now becomes that $a^q - a$ admits divided powers, as $(p)$ always admits divided powers.
One can use the natural surjection $D \to R/(X-a)$ to see that an element $a' \in A$
admits divided powers if and only if its image in $D$ admits divided powers.
Therefore, the condition becomes that
$a^q - a \in A$ should admit divided powers for all $a \in A$.
The above implies that in $A/p$ we have $(x^q - x)^p = 0$ for all $x \in A/p$,
since $A/p$ is assumed to be reduced, this is equivalent to all of its elements satisfying $x^q = x$.
Now we use \Cref{F_q lemma} to conclude that $A/p$ is actually a subalgebra of $\mathbb{F}_q$,
hence $A$ must be the Witt ring of a perfect subfield inside $\mathbb{F}_q$.
\end{proof}

\begin{remark}
Our argument excludes the existence of $q$-Frobenius if there is an element $a \in A$
such that $a^q - a$ does not admit divided powers.
For instance, if $A/p$ has a transcendental element over $\mathbb{F}_p$, then there is no functorial
endomorphism except for identity, as claimed in \cite[Remark 3.14.(3)]{LL20}.
It remains unclear to us, for instance, if the $p$-Frobenius can exist when $A = \mathbb{Z}_p[\sqrt{p}]$.
\end{remark}

Next we turn to Case II, which concerns the (base-changed) de Rham cohomology theory on algebras over a perfect ring of char.~$p$.
Once again the quasi-syntomic descent approach helps us prove the following statement (\textit{c.f.}~\cref{crystalline}):
\begin{proposition}
\label{Prop 3.7}
Let us assume that either
\begin{enumerate}
\item $A$ is an $\mathbb{F}_p$-algebra, $B$ is an $A$-algebra,
and consider the cohomology theory $\mathrm{dR} \otimes_A B$; or
\item $A = k$ is a perfect $\mathbb{F}_p$-algebra, $B$ is a $W_m(k)$-algebra,
and consider the cohomology theory $\mathrm{dR}_{m,1} \otimes_{W_m(k)} B$.
\end{enumerate}
Then the endomorphism monoids of the cohomology theory is a submonoid of $\mathbb{N}(\mathrm{Spec}(B))$,
where $\mathbb{N}$ stands for the constant monoid scheme of natural numbers with $1$ corresponding to the Frobenius.
\end{proposition}

\begin{proof}
We largely follow the strategy from the proof of \cite[Theorem 3.13]{LL20}.
Let us temporarily denote the cohomology theory by $\mathcal{F}$.

Note that in both cases, the functor $\mathcal{F}$ defines a quasi-syntomic sheaf on $\mathrm{qSyn}_A$:
for case (1) this is \Cref{p-completed base changed dR is a qsyn sheaf},\footnote{Since $A$ is $p$-torsion, we can drop the $p$-completion of the tensor product} and for case (2) this is \Cref{base change dR is a qsyn sheaf}.
Therefore, we can restrict ourselves to the category of QRSP $A$-algebras to compute the endomorphism monoid.

Next we reduce to one particular QRSP $A$-algebra: $R = A[X^{1/p^\infty}]/(X)$.
To make the reduction, apply the trick in proof of \cite[Proposition 7.10]{BS19} or \cite[Theorem 3.13]{LL20}
to see that for any QRSP $A$-algebra $S$ there exists an explicit QRSP $A$-algebra $S' = A[X_i^{1/p^\infty}; i \in I]/(f_j; j \in J)$,
where $f_j$ is an ind-regular sequence in $A[X_i^{1/p^\infty}; i \in I]$, together with a surjection
$R' \to R$ inducing a surjection of their values of the cohomology theory.
Hence for any functorial endomorphism,
its effect on $\mathcal{F}(S)$ is determined by that on $\mathcal{F}(S')$.
Finally for each $j \in J$, there exists a map $R \to S'$ sending $X^{\ell/p^n} \mapsto (f_j^{1/p^n})^{\ell}$,
the image of $\mathcal{F}(R)$ under these maps generates $\mathcal{F}(S')$;
therefore, the effect of a functorial endomorphism is determined by its effect on $\mathcal{F}(R)$.

Lastly we need to understand the effect of a potential functorial endomorphism $f$ on 
$D:= \mathcal{F}(R) = D_{(x)}(B[x^{1/p^\infty}])$, the divided power envelope of $(x)$ in $B[x^{1/p^\infty}]$.
From the proof of \cite[Theorem 3.13]{LL20} (last 4 paragraphs), we see that there is a finite disconnection
of $\mathrm{Spec}(B)$, such that on the $j$-th component we have $f(x^\ell) = x^{\ell \cdot p^{n_j}}$
for some natural number $n_j$.
Arguing component-wise, we may assume without loss of generality that $f(x^\ell) = x^{\ell \cdot p^{N}}$
for some natural number $N$; we need to show that this extends uniquely (assuming functoriality) to the whole $D$.
The algebra $D$ admits a natural grading, inspecting with the functoriality given by the map
$R \to R \otimes_A A[t^{1/p^\infty}]$ sending $x^{\ell}$ to $x^{\ell} \otimes t^{\ell}$ shows that $f$ must multiply the degree
by $p^N$.
Now we claim for every $n \in \mathbb{N}$ the effect of $f$
on \{degree $< p^{n+1}$ parts of $D$\} is determined by the effect of $f$ on \{degree $< p^{n}$ parts of $D$\},
which will finish the proof.
To that end, notice that the degree $< p^{n+1}$ parts is generated by $\gamma_{p^{n+1}}(x)$ and
the degree $< p^n$ parts.
Finally, we look at the map $A[x^{1/p^\infty}]/(X) \to A[y^{1/p^\infty}, z^{1/p^\infty}]/(y, z)$ 
given by $x^{\ell/p^i} \mapsto (y^{1/p^i} + z^{1/p^i})^\ell$.
By comparing the coefficients of 
$\gamma_{p^{n+N}}(y) \cdot \gamma_{p^{n+N}(p-1)}(z)$ of the equation obtained from functoriality,
one sees that the effect of $f$ on $\gamma_{p^{n+1}}(x)$ is pinned down by its effect on $\gamma_{p^{n}}(x)$
and $\gamma_{p^{n}(p-1)}(x)$.
\end{proof}

To illustrate the last sentence of the above proof, let us take $n = 0$ and see how to pin down the
effect of $f$ on $\gamma_p(x)$. The functoriality gives us a commutative diagram
\[
\xymatrix{
D \ar[r]^{f} \ar[d]_{x \mapsto (y+z)} & D \ar[d]^{x \mapsto (y+z)} \\
D \otimes_B D \ar[r]^{f \otimes f} & D \otimes_B D.
}
\]
Tracing through commutativity for the element $\gamma_p(x)$ we get that,
if $f(\gamma_p(x)) = c \cdot \gamma_{p^{N+1}}(x)$
then we have
\[
c \cdot \gamma_{p^{N+1}}(y) + \sum_{1 \leq j \leq p-1} \frac{1}{j!(p-j)!} y^{p^N} \cdot z^{p^N(p-1)} +
c \cdot \gamma_{p^{N+1}}(z) =
c \cdot \sum_{i + j = p^{N+1}} \gamma_i(y) \gamma_j(z).
\]
Therefore, we get $y^{p^N} \cdot \frac{z^{p^N(p-1)}}{(p-1)!} = c \cdot \gamma_{p^N}(y) \cdot \gamma_{p^N(p-1)}(z)$
in $D \otimes_B D$, which clearly pins down 
\[
c = \frac{(p^N)! \cdot (p^N(p-1))!}{(p-1)!}.
\]

Similar to \Cref{integral endomorphism prop}, if we make reducedness assumption on $B/p$ then we can
further decide which powers of Frobenius can appear depending on the size of $B/p$.
In \Cref{crystalline}, 
using the stacky approach, we will say precisely what powers of Frobenius
is allowed in terms of the map $k \to B^{\flat}$ for the case (2) in \Cref{Prop 3.7},
see \Cref{Compute Frob-k}.

\section{Endomorphisms of de Rham cohomology II}
\label{endo}

In this section we use a stacky approach to calculate endomorphism of de Rham and crystalline cohomology functors
in situations where it seems difficult to use only quasi-syntomic descent methods to obtain the answer.

\subsection{Unwinding equivalence}\label{sec4}

We fix two integers $n,m \ge 1$ and a perfect ring $k$ as before. The goal of this section is to study endomorphisms of the functor 
$$\mathrm{dR}_{m,n}: \w{A}\w{Rings}_{W_n(k)} \to \w{CAlg}(D(W_m(k))).$$ First, we will formulate this as a moduli problem. Let $S$ be a discrete test $W_m(k)$-algebra. We can define a functor $\mathrm{End}_{m,n}$ by $$\mathrm{End}_{m,n} (S):= \mathrm{End}(\mathrm{dR}_{m,n} \otimes_{W_m(k)} S).$$ This defines a functor $\mathrm{End}_{m,n}$ from $W_m(k)$-algebras to spaces which \textit{a priori} is a prestack. 
Let us study the base-changed crystalline cohomology theory, similar to \Cref{p-completed base changed dR is a qsyn sheaf} we have
the following 

\begin{proposition}
\label{base change dR is a qsyn sheaf}
The functor $\mathrm{dR}_{m,n} \otimes_{W_m(k)} S$, when restricted to $\mathrm{qSyn}_{W_n(k)}$, defines a quasi-syntomic sheaf.
\end{proposition}{}

\begin{proof}
Let us denote the derived crystalline cohomology functor relative to $W$ by $\mathrm{dR}_{\infty, n}$,
then we have $\mathrm{dR}_{m,n} \otimes_{W_m(k)} S \simeq \mathrm{dR}_{\infty, n} \otimes_{W(k)} S$.
Using the previous description and the fact that $W(k)$ is $p$-torsion free, to check quasi-syntomic sheaf
property it suffices to derived modulo $p$.
Since $(\mathrm{dR}_{\infty, n} \otimes_{W(k)} S)/^{L} p \simeq \mathrm{dR}_{1,n} \otimes_k (S/^{L} p)$,
we may reduce to the case where $m = n = 1$ and $S$ is a $1$-truncated animated $k$-algebra.
The proof of \Cref{p-completed base changed dR is a qsyn sheaf} works verbatim in this setting as well.
\end{proof}

\begin{lemma}
\label{dRmn discrete endo space}
The space of endomorphisms $\mathrm{End}(\mathrm{dR}_{m,n} \otimes_{W_m(k)} S)$ is discrete.
\end{lemma}

\begin{proof}
Similar to the proof of \Cref{discrete endo space}:
since $\mathrm{dR}_{m,n} \otimes_{W_m(k)} S$ defines a quasi-syntomic sheaf by \Cref{base change dR is a qsyn sheaf},
the claim follows from the fact that for a large quasi-syntomic $W_n(k)$-algebra $R$,
the value $\left(\mathrm{dR}_{m,n} \otimes_{W_m(k)} S\right)(R) = \mathrm{dR}_{m,n}(R) \otimes_{W_m(k)} S$ is a discrete algebra.
\end{proof}

On the other hand, let us consider the stack $\mathbb{A}^{1,\mathrm{dR}}$ which will always be viewed as a $W_n(k)$-algebra stack over $W_m(k)$ in this section. We define the following prestack capturing the endomorphisms of this stack along with the extra algebra structure.
\begin{notation}
For a test $W_m(k)$-algebra $S$, let us use $\mathscr{S}_{m,n}(S)$ to denote the space (groupoid) of endomorphisms of the stack $\mathbb{A}^{1, \mathrm{dR}}_{(S,n)}:=\mathbb{A}^{1, \mathrm{dR}}_{m,n} \times_{\w{Spec}\, W_m(k)} \w{Spec}\, S$ as a $W_n(k)$-algebra stack over $\w{Spec}\, S.$ 
\end{notation}{}

\begin{proposition}[Unwinding equivalence]\label{unwindequiv}The unwinding functor induces an isomorphism of prestacks 
$$\mathrm{Un}: \mathscr{S}_{m,n} \simeq \mathrm{End}_{m,n}. $$

\end{proposition}{}

\begin{proof}Unwinding provides a map from the left hand side to the right hand side. To show that this map is an isomorphism, let us fix a test $W_m(k)$-algebra $S.$ The $W_n(k)$-algebra stack $\mathbb{A}^{1, \mathrm{dR}}_S$ by definition
is an object of $\w{Fun}_{\times}(\w{Poly}_{W_n(k)}^{\w{op}}, \w{Stacks}_S).$ Since $\mathbb{A}^{1, \mathrm{dR}}_S$ is an affine stack (\cref{affinestack}) and the category $\w{AffStacks}_S$ is a full subcategory of $\w{Stacks}_S$ which is closed under small limits, we note that $\mathbb{A}^{1,\mathrm{dR}}_S$ is classified by an object of the full subcategory $\w{Fun}_{\times}(\w{Poly}_{W_n(k)}^{\w{op}}, \w{AffStacks}_S).$ By \cref{affcosimp}, the global section functor induces an equivalence of $\infty$-categories $\w{AffStacks}_S \simeq \w{coSCR}_S^{\w{op}},$ where the latter denotes the $\infty$-category of cosimplicial $S$-algebras. Therefore, $\mathbb{A}^{1,\mathrm{dR}}_S$ can be equivalently viewed as an object of $\w{Fun}(\w{Poly}_{W_n(k)}, \w{coSCR}_S)$. 
Therefore, endomorphisms of $\mathbb{A}^{1,\mathrm{dR}}_S$ as a $W_n(k)$-algebra stack can be computed as endomorphisms of the classifying object which we may call $G$ inside the category $\w{Fun}(\w{Poly}_{W_n(k)}, \w{coSCR}_S)$.

Now we look at the $S$-valued points of $\mathrm{End}_{m,n}.$ By property of left Kan extensions, this is given by endomorphisms of $\mathrm{dR}_{m,n} \otimes_{W_m(k)} S$ as a functor from $\w{Poly}_{W_n(k)} \to \w{CAlg}(D(S)).$ We can also left Kan extend along the inclusion $\w{Poly}_{W_n(k)} \to \mathrm{qSyn}_{W_n(k)}$ and equivalently consider endomorphisms of the functor $H:\mathrm{qSyn}_{W_n(k)} \to \w{CAlg}(D(S)).$ By \Cref{base change dR is a qsyn sheaf} we see that $H$ is a quasi-syntomic sheaf.

 A basis for the quasisyntomic topology on $\mathrm{qSyn}_{W_n(k)}$ is given by flat algebras over $W_n(k)$ whose reduction modulo $p$ is a QRSP algebra over $k.$ The category of such algebras will be denoted as $\w{QRSP}_{W_n(k)}$. On such algebras, the functor $H$ takes values in discrete rings. By properties of right Kan extension, we obtain that the functor $H$ has a canonical enrichment as a functor $H: \mathrm{qSyn}_{W_n(k)} \to \w{coSCR}_S$ and endomorphisms can also be calculated in the category $\w{Fun}(\mathrm{qSyn}_{W_n(k)}, \w{coSCR}_S).$ 
 By \Cref{identify endomorphisms}, we see that restricting along $\w{Poly}_{W_n(k)} \to \mathrm{qSyn}_{W_n(k)}$ now realizes $G$ as the canonical enrichment of 
 $\mathrm{dR}_{m,n} \otimes_{W_m(k)} S$.
 By property of left Kan extension, the endomorphisms of $H$ can also be computed as endomorphisms of $G$ in the category $\w{Fun}(\w{Poly}_{W_n(k)}, \w{coSCR}_S),$ 
 which finishes the proof.
\end{proof}{}

\begin{proposition}The functor $\mathscr{S}_{m,n} : \mathrm{Alg}_{W_m(k)} \to \mathcal{S}$ is an fpqc sheaf. In fact, it is a sheaf of sets. 

\end{proposition}{}
\begin{proof}
This follows from \Cref{dRmn discrete endo space} and the fact that $\mathbb{A}^{1,\w{dR}}$ is an fpqc stack.
\end{proof}{}

Before we proceed further, let us make the following definition. Let $m \ge 1$ be an arbitrary integer fixed as before. Then $\mathbb{G}_a^{\w{perf}}$ represents an fpqc sheaf of rings on the category of $W_m(k)$--algebras. 

\begin{definition}We define a sheaf $\underline{\mathrm{Frob}}_k: \w{Alg}_{W_m(k)} \to \w{Sets}$ to be the subsheaf of $\mathcal{H}om_{k \w{-}alg}(\mathbb{G}_a^{\w{perf}},
\mathbb{G}_a^{\w{perf}})$ such that if $B$ is a $W_m(k)$--algebra, then $\underline{\mathrm{Frob}}_k(B)$ is the set of $k$-algebra scheme maps $\mathbb{G}_{a,B}^{\w{perf}} \to \mathbb{G}_{a,B}^{\w{perf}}$ which is induced by an algebra map $B[x^{1/p^\infty}] \to B[x^{1/p^\infty}]$ such that it sends $x \to \sum_{i} b_i x^{p^i}$ where the sum ranges over a finite subset in $\mathbb{Z}_{\geq 0}$. The sheaf $\underline{\mathrm{Frob}}_k$ naturally has the structure of a commutative monoid.

\end{definition}{}

\begin{notation}For a $W_m(k)$--algebra $B,$ we will write the symbol $\w{Frob}^i$ to mean an element of $i \in \underline{\mathrm{Frob}}_k(B).$ We will also write $\w{Frob}^{i+j}$ to denote the composition of $\w{Frob}^i$ and $\w{Frob}^j.$ 

\end{notation}{}

\begin{remark}
\label{Compute Frob-k}
We note that $\underline{\mathrm{Frob}}_k$ is a subsheaf of the sheafification of the constant monoid $\underline{\mathbb{N}}.$ In fact, they are equal when $k = \mathbb{F}_p,$ but this is not always the case.
One can compute that, given a $W_m(k)$-algebra $B$, we have
\[
\underline{\mathrm{Frob}}_k(B) = \w{Hom}_k (\mathbb{G}_{a, B^{\flat}}, \mathbb{G}_{a, B^{\flat}}).
\]
In very concrete terms, the right hand side above
is the set of pairs $(\mathcal{P}, i)$, 
where $\mathcal{P}$ is a partition $B = \prod_{j \in J} B_j$ and $i = (i_j)$
is a function on $\mathrm{Spec}(B)$,
which is constant on each $\mathrm{Spec}(B_j)$ taking values in $\mathbb{N}$,
satisfying the condition that the map $W_m(k)^{\flat} = k \to B_j^{\flat}$
factors through a subfield of the finite field $\mathbb{F}_{p^{i_j}}$.

Consequently, one finds that when $k$ is a perfect \emph{field},
the sheaf $\underline{\mathrm{Frob}}_k$ is representable by either 
the constant monoid scheme $\mathbb{N}$ or the singleton $\{0\}$,
depending on whether $k$ is finite or not.
\end{remark}{}

\begin{proposition}\label{trivial1} There is an isomorphism of sheaf of monoids over $k$ $$\mathscr{S}_{1,1} \simeq \underline{\mathrm{Frob}}_k. $$ 
\end{proposition}{}
\begin{proof}
Let $k$ be a perfect ring. Let $B$ be an arbitrary $k$-algebra. By \cref{Compute Frob-k}, our goal is to show that $\mathrm{End}(\w{dR} \otimes^{L}_k B)$ is just given by $\w{Hom}_{k} (\mathbb{G}_{a,B}, \mathbb{G}_{a,B})$ where $\mathbb{G}_{a,B}$ is regarded as a $k$-algebra scheme over $B.$ For the proof, we will use another $k$-algebra which we denote as $\mathbb{G}_{a,B}^{\w{perf}}.$ More explicitly, $\mathbb{G}_a^{\w{perf}}$ is represented by the affine scheme $\w{Spec}\, B[x^{1/p^\infty}]$ (see \cref{gaperf}). Note that we have a natural injection of sets $i: \w{Hom}_{k} (\mathbb{G}_{a,B}, \mathbb{G}_{a,B}) \to \w{Hom}_{k} (\mathbb{G}_{a,B}^{\w{perf}}, \mathbb{G}_{a,B}^{\w{perf}}).$

Let us first construct a map
$\varphi: \mathrm{End}(\w{dR} \otimes^{L}_k B) \to \w{Hom}_{k} (\mathbb{G}_{a,B}^{\w{perf}}, \mathbb{G}_{a,B}^{\w{perf}}).$ We note that $\w{dR}$ restricts to a functor on the full subcategory of $k$-algebras which we denote as $\w{Poly}_{\w{perf}}/_k$ which consists of perfections of finite type polynomial algebras 
over $k.$ If $R \in \w{Poly}_{\w{perf}}/_k, $ then $\w{dR}_{R/k} \otimes_k B \simeq R \otimes_k B,$ which defines a functor from $\w{Poly}_{\w{perf}}/_k \to \w{Alg}/_B$ sending $R \mapsto R \otimes_k B.$ This basically classifies perfect $k$-algebra ring schemes over $\w{Spec}\, B$ and any endomorphism of $\w{dR} \otimes^{L}_k B$ induces an endomorphism of this perfect $k$-algebra ring scheme over $\w{Spec}\,B,$ which is just given by $\mathbb{G}_{a,B}^{\w{perf}}.$ This constructs the required map $\varphi.$

We know that any element in 
$\mathrm{End}(\w{dR} \otimes^{L}_k B)$ is uniquely determined by a map $f$ of $\mathbb{A}^{1, \w{dR}}$ as a $k$-algebra stack over $\w{Spec}\, B.$ We also note that there is a natural map $\mathbb{G}_a^{\w{perf}} \to \mathbb{A}^{1, \w{dR}}$ of $k$-algebra stacks over $\w{Spec}\, B$ (from now on, we will omit the $B$ from the subscript to ease our notation). By functoriality of $S \mapsto S_{\w{perf}}$ and the fact that this perfection construction commutes with colimits, it follows that the map $f$ lifts to give a map as below.
\begin{center}
    \begin{tikzcd}
\mathbb{G}_a^{\text{perf}} \arrow[d] \arrow[r, "\hat{f}"] & \mathbb{G}_a^{\text{perf}} \arrow[d] \\
{\mathbb A ^{1, \text{dR}}} \arrow[r, "f"]                & {\mathbb A ^{1, \text{dR}}}         
\end{tikzcd}
\end{center}{}
Let $u: \mathbb{G}_a^{\w{perf}} \to \mathbb{G}_a$ denote the natural map of $k$-algebra schemes. Then fibre of the map $\mathbb{G}_a^{\w{perf}} \to \mathbb{A}^{1, \w{dR}}$ identifies with $u^* W[F]$ (\textit{c.f.}~\cite[Proposition 2.2.6]{Mon21}). Therefore, $f$ is given by a map of the quasi-ideal in $\mathbb{G}_a^{\w{perf}}$ given by $u^* W[F] \to \mathbb{G}_a^{\w{perf}},$ which is of the form of a commutative diagram as below.
\begin{center}
    \begin{tikzcd}
{u^* W[F]} \arrow[r] \arrow[d, "t"'] & \mathbb{G}_a^{\text{perf}} \arrow[d, "\varphi(f)"] \\
{u^* W[F]} \arrow[r]                 & \mathbb{G}_a^{\text{perf}}                        
\end{tikzcd}
\end{center}{}
In the above, $t$ is required to be a $\mathbb{G}_a^{\w{perf}}$-module map once the target is given the appropriate $\mathbb{G}_a^{\w{perf}}$-module structure via restricting scalars along $\varphi(f).$ Now, inspecting the above diagram at the level of global sections yields that the map $\varphi$ must factor through $i,$ i.e., $\varphi(f)$ must be induced by an element of $s \in \w{Hom}_k(\mathbb{G}_a, \mathbb{G}_a).$ From this, it follows that the previous commutative diagram is uniquely determined by a commutative diagram as below.
\begin{center}
    \begin{tikzcd}
{W[F]} \arrow[d, "t'"'] \arrow[r] & \mathbb{G}_a \arrow[d, "s"] \\
{W[F]} \arrow[r]                  & \mathbb{G}_a               
\end{tikzcd}
\end{center}{}
In the above $t$ is required to be a $\mathbb{G}_a$-module map once the target is given the appropriate $\mathbb{G}_a$-module structure via restricting scalars along $\varphi.$ In order to understand the map $t',$ we can therefore apply graded Cartier duality \cite[\S 2.4]{Mon21}. We note that $W[F]^* = \mathbb{G}_a$ and thus we get a map of graded group schemes $t'^*: \mathbb{G}_a \to \mathbb{G}_a,$ where the source group scheme $\mathbb{G}_a$ receives its grading via the $\mathbb{G}_a$-module structure induced by restriction of scalars along $s.$ By easy degree considerations, it follows that there exists a unique $\mathbb{G}_a$-module map $t'$ which fits into the above commutative diagram. Therefore, we obtain the natural bijection $\mathrm{End}(\w{dR} \otimes^{L}_k B) \simeq \w{Hom}_k(\mathbb{G}_{a,B}, \mathbb{G}_{a,B}),$ as desired.  
\end{proof}{}

\begin{proposition}\label{crystalline}For any $m \ge 1$, there is a natural isomorphism of sheaf of monoids over $W_m(k)$
$$\mathscr{S}_{m,1} \simeq \underline{\mathrm{Frob}}_k.$$ 

\end{proposition}{}

\begin{proof}
Let $B$ be a $W_m(k)$-algebra. There is a $k$-algebra scheme over $B$, which we denote as $\mathbb{G}_{a,B}^\w{perf}$, whose underlying affine scheme is $\w{Spec}\, B [x^{1/p^\infty}].$ As in the proof of \cref{trivial1}, one also obtains a map $\varphi: \mathrm{End}(\w{dR} \otimes_{W_m(k)} B) \to \w{Hom}_{k} (\mathbb{G}_{a,B}^{\w{perf}}, \mathbb{G}_{a,B}^{\w{perf}}).$ It follows from going modulo $p$ and applying \cref{trivial1} that $\varphi$ actually factors to give a map again denoted as $\varphi: \mathrm{End}(\w{dR} \otimes_{W_m(k)} B) \to \underline{\mathrm{Frob}}_k(B).$ We will argue that this map is a bijection.

By using the stack $\mathbb{A}^{1,\w{dR}}$ and the natural map $\mathbb{G}_a^{\w{perf}} \to \mathbb{A}^{1, \w{dR}},$ in a way similar to the proof of \cref{trivial1}, this amounts to the more concrete assertion that there is a unique map $t$ of quasi-ideals in $\mathbb{G}_a^{\w{perf}}$ as below.

\begin{center}
    \begin{tikzcd}
{\mathbb{G}_a^{\w{perf}, \sharp}} \arrow[r] \arrow[d, "t"'] & \mathbb{G}_a^{\text{perf}} \arrow[d, "\iota(x)"] \\
{\mathbb{G}_a^{\w{perf}, \sharp}} \arrow[r]                 & \mathbb{G}_a^{\text{perf}}                        
\end{tikzcd}
\end{center}{}
Here $x \in \underline{\mathrm{Frob}}_k(B)$ and $\iota: \underline{\mathrm{Frob}}_k(B) \to \w{Hom}_{k} (\mathbb{G}_{a,B}^{\w{perf}}, \mathbb{G}_{a,B}^{\w{perf}})$ denotes the natural inclusion.
Let us write $U$ for the coordinate ring of $\mathbb{G}_a^{\w{perf}, \sharp}.$ Then $U$ is a $\mathbb{N}[1/p]$-graded Hopf algebra over $B.$ It is also a free algebra over $B$, where all the homogeneous components are free of rank $1$ over $B.$ As a graded $B$-algebra, $U$ is generated by the basis elements in degree $p^i$ for $i \in \mathbb{Z}.$ It is enough to check that for a fixed $x \in \underline{\mathrm{Frob}}_k(B),$ there exists a unique map $t$ which gives a map of quasi-ideals as above. The existence is clear from definition of $\mathbb{G}_a^{\w{perf}, \sharp}$ (\cref{gaperfsharp}) by applying the divided power envelope construction. For the uniqueness we note that once $x$ is fixed, the above diagram forces the homogeneous elements of degree $p^i$ for $ i \le 0$ to be mapped uniquely. 
The rest follows from inspecting the comultiplication of $U$ and induction on $i$ (see last paragraph of the proof of \cref{Prop 3.7}
as well as the discussion after that proof).   
\end{proof}{}

\begin{remark}
It is possible to prove \cref{crystalline} by using the methods from \cite[\S 3.4]{Mon21}, which would essentially amount to proving a similar statement as above about the quasi-ideal $\mathbb{G}_a^{\w{perf}, \sharp} \to \mathbb{G}_a^{\w{perf}}.$ It is also possible to reduce to the same statement about quasi-ideals directly from \cref{perfectionlemma} by using the compatibility of the map induced on the animated ring $W(S)/^Lp$ via the Frobenius on $W(S)$ with the natural Frobenius operator on any animated $k$-algebra: this implies that any endomorphism of $\mathbb{A}^{1, \w{dR}}$ as a $k$-algebra stack lifts along the map $\mathbb{G}_a^{\w{perf}} \to \mathbb{A}^{1, \w{dR}}$ obtained by taking perfection. This lifting property fails for endomorphisms of $\mathbb{A}^{1, \w{dR}}$ as a $W_2(k)$-algebra stack, leading to extra endomorphisms as will be constructed in \cref{enough}.
\end{remark}{}

\subsection{Construction of endomorphisms}\label{enough}
This subsection describes the construction of ``enough'' endomorphisms of de Rham cohomology. Our strategy is to crucially exploit the Unwinding equivalence proven in \cref{unwindequiv} to pass to the world of ring stacks and do a small explicit construction there. We will begin fixing notations and making some definitions. Since we are interested in endomorphisms, we will ignore the Frobenius twist introduced in \cref{frobtwist}.

\begin{notation}In this section, we will work with a perfect ring $k$ of characteristic $p >0.$ 
We fix two integers $n,m \ge 1$.
We will use $W$ to denote the Witt ring scheme over the fixed base $W_m(k).$ 
Since $m,n$ are fixed, we will denote $\mathbb{A}^{1, \mathrm{dR}}_{(m,n)}$ 
simply by $\mathbb{A}^{1, \mathrm{dR}}$ when no confusion is likely to occur.

\end{notation}{}

\begin{definition}We will let $W[p]$ denote the group scheme underlying the kernel of multiplication by $p$ map on $W.$
\end{definition}{}

\begin{definition}
\label{notation of aut gp}
We will let $(1 + W[p])^\times$ denote the monoid scheme underlying $x \in W$ satisfying $px = p.$ The multiplication on this monoid scheme is given by simply using the multiplication underlying the ring scheme structure on $W.$

\end{definition}{}

\begin{proposition}
\label{the monoid scheme is group}
Let $B$ be a $p$-nilpotent ring, then
the monoid scheme $(1+ W[p])^\times$ over $\mathrm{Spec}(B)$ is a group scheme.

\end{proposition}{}

\begin{proof} 
This amounts to saying that for any ring $S$ with $p^m = 0$ in $S$
for some $m$,
if $x \in W(S)$ satisfies $px = p$ then $x$ must be a unit in the ring $W(S).$ 
Recall that we have a short exact sequence
\[
0 \to W(p \cdot S) \to W(S) \to W(S/p) \to 0,
\]
where $W(p \cdot S)$ denotes the Witt ring associated with the ideal 
(viewed as a non-unital ring) $p \cdot S$.
Since $p^m = 0$ in $S$, we know that the ideal
$W(p \cdot S)$ is nilpotent.
Therefore, it suffices to show that the image of $x$ in $W(S/p)$ is a unit,
hence we have reduced to the case where $S$ is of characteristic $p$.
Since $p = V(1)$ in this case, the condition of $x$ reads
$V(F(x)) = x \cdot V(1) = V(1)$.
Injectivity of $V$ shows that $F(x) = 1$, which implies that $x$ is a unit.
\end{proof}{}

\begin{construction}
\label{construct aut}
Now we will begin our construction of endomorphisms of $\mathbb{A}^{1, \mathrm{dR}}$ as a $W_n(k)$-algebra stack (over the base $W_m (k)$ which is fixed for this section) when $n \ge 2$. Since \cref{a1dR} constructs the above stack as cone of the quasi-ideal $d: W \xrightarrow[]{\times p} W,$ we will explicitly construct maps at the quasi-ideal level, which can be done purely $1$-categorically. We note that there is a natural structure map $(W(k) \xrightarrow[]{\times p} W(k) ) \to  (W \xrightarrow[]{\times p} W)$ of quasi-ideals, which describes the structure of $(W \xrightarrow[]{\times p} W)$ as a quasi-ideal over $k.$ In the language of quasi-ideals, the natural map $W_n(k) \to k$ which gives that $W_n(k)$ structure map of $(W(k) \xrightarrow[]{\times p} W(k) )$ viewed as a map of quasi-ideals as described below.

\begin{center}
    \begin{tikzcd}
W \arrow[r, "\times p"]                           & W              \\
W(k) \arrow[r, "\times p"] \arrow[u]              & W(k) \arrow[u] \\
W(k) \arrow[r, "p^n"] \arrow[u, "\times p^{n-1}"] & W(k) \arrow[u]
\end{tikzcd}
\end{center}{}

We will construct maps of the quasi-ideal $d: W \xrightarrow[]{\times p} W$ over the structure map $W(k) \xrightarrow[]{\times p^n} W(k)$ as described above. Let $F$ be a homomorphism of the $W(k)$-algebra scheme $W.$ A quasi-ideal map from $d: W \xrightarrow[]{\times p} W$ to itself can be defined by giving a $W$-linear map $u: W \to F_* W$ which makes the diagram below commutative. 

\begin{center}
    \begin{tikzcd}
W \arrow[d, "u"] \arrow[r, "\times p"] & W \arrow[d, "F"] \\
W \arrow[r, "\times p"]                & W               
\end{tikzcd}
\end{center}{}

However, we need to make sure that such a map respects the additional structure of being a map of quasi-ideals over $W(k) \xrightarrow[]{\times p^n} W(k),$ i.e., the following diagram needs to commute.

\begin{center}
    \begin{tikzcd}
                                                                                             &                                          & W \arrow[d, "u"] \arrow[r, "\times p"] & W \arrow[d, "F"] \\
                                                                                             &                                          & W \arrow[r, "\times p"]                & W                \\
W(k)  \arrow[rruu, "p^{n-1}", bend left] \arrow[r, "\times p^{n}"'] \arrow[rru, "\times p^{n-1}"] & W(k) \arrow[rruu, bend left] \arrow[rru] &                                        &                 
\end{tikzcd}
\end{center}{}
As one checks, for any $n \ge 2$, the only condition this imposes is that $p u(1) = p.$ This provides the following map that we wanted to construct.

$$(1+ W[p])^\times \cdot F \to \mathscr{S}_{m,n} $$
Further, for any $n\ge 2,$ the above map is clearly an injection by construction. We point out that it is possible to do such a construction for every $W(k)$-algebra map $F$ of the ring scheme $W.$ Let $S$ be a $W_m(k)$-algebra. Then the element of $\mathscr{S}_{m,n}(B)$ constructed above will be denoted by $u \cdot F$ where $u$ would be understood to be an element $u(1) \in (1+ W[p])^\times (B).$
By construction, we see that the composition $(u, F') \circ (v,F)$ is equal to $(u F'(v), F'F).$
\end{construction}{}

\begin{remark}Note that in the above picture if we let $n=1,$ then $u(1)$ is \textit{forced by the diagram} to be equal to $1,$ and one does not get the extra endomorphisms that was constructed above for $n \ge 2.$

\end{remark}{}

\begin{proposition}
\label{automorphism proposition}
Let $(1+ W[p])^\times$ denote group scheme as above. There is an injection of (sheaves) 
$\coprod_{i \in \underline{\mathrm{Frob}}_k} (1+ W[p])^\times \cdot \mathrm{Frob}^i \to \mathrm{End}_{m,n}$
when $n \geq 2$.
\end{proposition}{}

\begin{proof}This follows from \cref{unwindequiv} and \cref{construct aut}.
\end{proof}{}

\begin{remark}
\label{need}
Let $B$ be a $W_m(k)$-algebra;
we construct two natural maps
\[
\underline{\mathrm{Frob}}_k(B) \to \mathrm{End}_{W(k)}(W^{(1)}_B) \to \underline{\mathrm{Frob}}_k(B).
\]

The first arrow follows from the explicit description given in \Cref{Compute Frob-k},
and we simply send powers of the Frobenius to powers of Frobenius on the Witt ring scheme.
To exhibit the second arrow, note that any element in $\mathrm{End}_{W(k)}(W^{(1)}_B)$
induces an element in $\mathrm{End}_k([W^{(1)}_B / p]) \simeq \mathrm{End}_k(\mathbb{A}^{1, \mathrm{dR}}_B)$
which is equivalent to $\underline{\mathrm{Frob}}_k(B)$ by \Cref{crystalline}.
One easily checks that the composition of the two maps gives identity on $\underline{\mathrm{Frob}}_k(B).$
\end{remark}

\subsection{Calculation of the endomorphism monoid}

Throughout this subsection, we will fix $k$ to be a perfect ring as before. In this subsection, we will show that we have found all the endomorphisms of $\w{dR}_{m,n}$; more precisely, the injection in \cref{automorphism proposition} is an isomorphism.

We need some preparations, starting with understanding the homotopy sheaves associated with
$\mathbb{A}^{1, \mathrm{dR}}$. Since $\mathbb{A}^{1, \mathrm{dR}}$ is a $1$-stack, we only need to understand $\pi_0$ and $\pi_1.$ Once again, we remind the readers that since we are interested in endomorphisms, we will ignore the Frobenius twist introduced in \cref{frobtwist}.

\begin{proposition}
For a test algebra $S,$
\label{cohomology of A1dR}
\begin{enumerate}
\item $\mathbb{A}^{1, \mathrm{dR}}(S) = W(S)/^{L} p$, where $W(S)/^{L} p$ denotes the animated ring obtained by quotienting $W(S)$ by $p.$ We note that the object in the category of animated modules underlying $W(S)/^{L} p$ can be simply described as $\mathrm{Cofib}(W(S) \xrightarrow[]{\times p} W(S)).$
\item The sheaf $\pi_1(\mathbb{A}^{1, \mathrm{dR}})$ is representable by $W[p]$,
the ideal scheme of $p$-torsion in the ring scheme $W$.
\item Over a characteristic $p$ base, the sheaf $\pi_0(\mathbb{A}^{1, \mathrm{dR}})$
is representable by $\mathbb{G}_a,$ where the induced map $W \to \mathbb{G}_a$
is given by the natural projection to $0$-th Witt coordinate.
\end{enumerate}
\end{proposition}

\begin{proof}
\,

\begin{enumerate}
    \item By definition, we need to prove that the presheaf $P(S) := W(S)/^Lp$ is already an fpqc sheaf of animated rings. It is enough to show that $P(S):= \mathrm{Cofib}(W(S) \xrightarrow[]{\times p} W(S))$ is a sheaf of animated modules. By noting that $\mathrm{Cofib}(W(S) \xrightarrow[]{\times p} W(S)) = \mathrm{fib}(W(S)[1] \xrightarrow[]{\times p} W(S)[1]),$ we see that it is enough to prove that the functor $Q(S) := W(S)[1]$ is a sheaf of connective animated modules. For this, we only need to show that $H^1_{fpqc}(\mathrm{Spec}\, S, W) = 0.$

To this end, we note that the sheaf $W = \varprojlim_{n} W_n$. 
By \cite[Example 3.1.7 and Proposition 3.1.10]{proetale} and the fact that $F$ on $W$
is faithfully flat, it follows that $W= R \varprojlim_{n} W_n$. Thus $R\Gamma_{fpqc}(\w{Spec}\, S, W) = R\varprojlim_{n} R\Gamma_{fpqc} (\w{Spec}\, S, W_n).$
Now one notes that $W_n$ has a finite filtration with the graded pieces being equal to $\mathbb{G}_a.$ Thus $R\Gamma_{fpqc}(\w{Spec}\, S, W) = R\varprojlim_{n} R\Gamma_{fpqc} (\w{Spec}\, S, W_n) = R\varprojlim_{n} \Gamma(\w{Spec}\, S, W_n)= R\varprojlim_{n} W_n(S) = \varprojlim_{n} W_n(S)= W(S).$ In particular, $H^1_{fpqc}(\mathrm{Spec}\, S, W) = 0,$ as desired.

\item Follows from (1).

\item In the Witt ring of a characteristic $p$ ring, we have $p = VF$.
Therefore, the conclusion follows since $F: W \to W$ is a fpqc surjection. \vspace{-5.5mm}\end{enumerate}
\end{proof}

In general, $\pi_0 (\mathbb{A}^{1, \w{dR}})$ is given by the sheaf of discrete $k$-algebras $W/p.$ However, if the base is not of characteristic $p,$ this sheaf stops being representable as noted below. Nevertheless, \cref{image is Frobenius powers} will helps us extract the necessary information from $\pi_0 (\mathbb{A}^{1, \w{dR}})$ relevant to us.

\begin{proposition}
Let $B$ be a ring such that $p \not\in (p^2)$,
let $S = \mathrm{Spec}(B)$.
The sheaf $$
\mathcal{F} \coloneqq \pi_0(\mathrm{Cone}(\mathbb{G}_{a,S}^\sharp \to \mathbb{G}_{a,S})) \simeq
\pi_0(\mathbb{A}^{1, \mathrm{dR}}) $$
is not representable by an algebraic space over $S$. 
\end{proposition} 
\begin{proof} The isomorphism follows from \Cref{two models}. Since both $\mathbb{G}_a$ and $\mathbb{G}_a^{\sharp}$ are affine schemes,
the hypothetical representing algebraic space would be quasi-compact and quasi-separated.
Below we show that there cannot be such a qcqs algebraic space.

It suffices to prove the statement for $B/p^2$ hence we may assume that $p^2 = 0$ in $B$.
Since the restriction of our sheaf to $B/p$-algebras
is represented by the affine scheme $\mathbb{G}_{a, B/p}$,
using \cite[\href{https://stacks.math.columbia.edu/tag/07V6}{Tag 07V6}]{stacks-project}
we see that the sheaf would in fact be represented by an affine scheme over $S$.
Let us denote its ring of function by $R$.
The natural map $\mathbb{G}_{a,S} \to \mathrm{Spec}(R)$ induces a map $R \to B[t]$.
Reducing the ring map modulo $p$, we see that the image is $B/p[t^p]$. This implies that an element of the form $t^p + p \cdot g$ must be in the image.
On the other hand we claim the image of the ring map itself is contained
in $\{f \in B[t] \mid f'(t) = 0\}$.
Indeed, the two compositions
$$\xymatrix{
\mathrm{Spec}(B[t,\epsilon]/\epsilon^2) 
\ar@/^/[r]^-{t \mapsto t} \ar@/_/[r]_-{t \mapsto t+\epsilon}
& \mathrm{Spec}(B[t]) \ar[r] & \mathcal{F}
}$$
yields the same map as $\epsilon \in B[t,\epsilon]/\epsilon^2$ admits divided powers.
This shows that the image of $R \to B[t]$ must be contained 
in the equalizer of the two maps $B[t] \rightrightarrows B[t,\epsilon]/\epsilon^2$.
The identification of this equalizer with those polynomials whose derivative
is zero follows from Taylor expansion.
Lastly to get a contradiction, one observes that if we let $f = t^p + p \cdot g$, 
then $f' \not= 0$ as $p \not\in (p^2)$; however, we had previously argued that $t^p + p \cdot g$ must be in the image.
\end{proof}

\begin{lemma}
\label{possible Fp structure on A1dR}
Let $B$ be a $W(k)$-algebra.
We have $W(B)[p] \simeq \mathrm{Hom}_{W(k)}(k, \mathbb{A}^{1, \mathrm{dR}})$,
where the right hand side denotes the space of maps as $W(k)$-algebra stacks over $B.$ 
Given $\beta \in W(B)[p]$ the corresponding homomorphism of sheaves is modeled by
\[
\xymatrix{
W \ar[r]^{\times p} & W \\
W(k) \ar[r]^{\times p} \ar[u]^{\times (1 + \beta)} & W(k) \ar[u].
}
\]
\end{lemma}

Here, the constant sheaf of $W(k)$-algebras given by $k$ is viewed as a $W(k)$-algebra stack over $B.$

\begin{proof}
Since $\mathbb{A}^{1, \w{dR}}$ is $1$-truncated, by \cref{square-zero} the right hand side is classified by
$\mathrm{Hom}_k(\mathbb{L}_{k/W(k)}, \pi_1(\mathbb{A}^{1, \mathrm{dR}})[1])$ $= W(B)[p]$.
Here in this identification we have used \Cref{cohomology of A1dR}.(2).
One checks easily that the maps we constructed in the last sentence exactly
corresponds to $\beta$ under the above identification, hence finishing the proof.
\end{proof}

Our last preparation is to understand those algebra homomorphisms in
$\mathrm{End}_{k}(\pi_0(\mathbb{A}^{1, \mathrm{dR}}))
$ which can be lifted to a $W_n(k)$-algebra homomorphism of $\mathbb{A}^{1, \mathrm{dR}}$.
It turns out that liftability as a $W_n(k)$-algebra stack for $n>1$ \textit{automatically} guarantees liftability as
a $k$-algebra stack as noted below.

\begin{lemma}
\label{image is Frobenius powers}
Let $B$ be a $W_m(k)$-algebra and let us consider $\mathbb{A}^{1, \mathrm{dR}}$ as a $k$-algebra stack over $B$.
The two natural maps $
\mathrm{End}_{W(k)}(\mathbb{A}^{1, \mathrm{dR}}) \to 
\mathrm{End}_{k}(\pi_0(\mathbb{A}^{1, \mathrm{dR}}))
$
and
$\mathrm{End}_{k}(\mathbb{A}^{1, \mathrm{dR}}) \to \mathrm{End}_{k}(\pi_0(\mathbb{A}^{1, \mathrm{dR}}))$ have the same image.
In particular, by \cref{crystalline},
we know that the image is naturally in bijection with the monoid $\underline{\mathrm{Frob}}_k(B)$.
\end{lemma}

\begin{proof}
The image of the first map clearly contains the image of second map.
Now given $f \in \mathrm{End}_{W(k)}(\mathbb{A}^{1, \mathrm{dR}})$,
by composing with the natural map $\iota:k \to \mathbb{A}^{1, \mathrm{dR}}$,
we get a natural map $f \circ \iota \colon k \to \mathbb{A}^{1, \mathrm{dR}}$
of $W(k)$-algebra stacks.
In \Cref{possible Fp structure on A1dR} we see that $f \circ \iota$ must be classified by some element $1+\beta \in 1+W(B)[p]$. By \Cref{the monoid scheme is group}, we can find an inverse $(1 + \beta)^{-1} \in (1+ W[p])^\times$,
one sees that the composition $(1 + \beta)^{-1} \circ f \circ \iota = \iota$.
Here we regard an element in $(1+ W[p])^\times$ as a $W(k)$-algebra automorphism of $\mathbb{A}^{1, \mathrm{dR}}$
by \Cref{construct aut}.
Since these elements in $(1+ W[p])^\times$ always induce identity on $\pi_0$, we see that
$(1 + \beta)^{-1} \circ f$ is a $k$-algebra automorphism lifting the same ring homomorphism
on $\pi_0(\mathbb{A}^{1, \mathrm{dR}})$ as $f$.
\end{proof}

\begin{theorem}
\label{representing endomorphism monoid}
Let $A = W_n(k)$ and suppose that $B$ is a $W_m(k)$-algebra, then we have
\[
\mathrm{End}(\mathrm{dR}_{m,n} \otimes_{W_m(k)} B) =
\begin{cases}
\coprod_{i \in \underline{\mathrm{Frob}}_k(B)} \mathrm{Frob}^i &  n = 1; \\

\vspace{1mm}
\coprod_{i \in \underline{\mathrm{Frob}}_k(B)}(1+ W[p])^\times(B) \cdot \mathrm{Frob}^i &  n \geq 2; \\

\end{cases}
\]
Here the multiplication law in the second case is given by
\[
(u \cdot \mathrm{Frob}^i) \cdot (v \cdot \mathrm{Frob}^j) =
u \cdot \mathrm{Frob}^i(v) \cdot \mathrm{Frob}^{i+j}.
\]
where $u, v \in (1+ W[p])^\times(B)$.

\end{theorem}

\begin{remark}
Note that these endomorphism spaces are all discrete, by \Cref{discrete endo space}.
The above theorem states that the map in \Cref{automorphism proposition}
is actually an isomorphism.
From the above calculation, we also conclude that the sheaf of endomorphism monoids is representable
if and only if the sheaf $\underline{\mathrm{Frob}}_k$ is representable.
This happens whenever $k$ is a perfect \emph{field},
in which case the representing scheme is a combination of the constant monoid scheme $\mathbb{N}$ and
the commutative group scheme $(1+ W[p])^\times$, depending on $k$ and $n$.
\end{remark}
\begin{proof}
When $n = 1$, this is proved in \cref{crystalline}. Below we will assume that $n \geq 2$.

Recall that in \Cref{unwindequiv} we have shown that the endomorphisms of our de Rham cohomology functor is the same as
the endomorphisms of the $W_n(k)$-algebra stack $\mathbb{A}^{1, \text{dR}}_B$ over $\mathrm{Spec}(B)$.
Since the category of $W_n(k)$-algebra stacks is equivalent to the category of 
sheaves of $W_n(k)$-animated algebras, see \Cref{switch}, we will compute the endomorphism of
$\mathbb{A}^{1, \text{dR}}_B$ viewed as a sheaf of $W_n(k)$-animated algebras on the fpqc site of ${\mathrm{Spec}(B)}$.

Composing with the map $\mathbb{A}^{1, \text{dR}} \to \pi_0(\mathbb{A}^{1, \text{dR}})$, we get a natural map
\[
\mathrm{Hom}_{W_n(k)}(\mathbb{A}^{1, \text{dR}}, \mathbb{A}^{1, \text{dR}}) \xrightarrow{f_n}
\mathrm{Hom}_{W_n(k)}(\mathbb{A}^{1, \text{dR}}, \pi_0(\mathbb{A}^{1, \text{dR}})) =
\mathrm{End}_{k}(\pi_0(\mathbb{A}^{1, \text{dR}})),
\]
here and below, all appearances of $\mathrm{Hom}$ refers to homomorphism of sheaves respecting the
designated structure marked by subscript.
By \Cref{image is Frobenius powers}, we see that
\[
\mathrm{Im}(f_n) = \underline{\mathrm{Frob}}_k(B)
\]

We need to understand the fibre of $f_n$.
Take an $i \in \underline{\mathrm{Frob}}_k(B)$, by \cref{square-zero},
the fibre of $f_n$ over $\mathrm{Frob}^i$ is a torsor under 
\[
\mathrm{Hom}_{\mathbb{A}^{1, \text{dR}}}(\mathbb{L}_{\mathbb{A}^{1, \text{dR}}/W_n(k)}, \pi_1(\mathbb{A}^{1, \text{dR}})[1]).
\]
Here the sheaf of $\mathbb{A}^{1, \text{dR}}$-module structure on the sheaf $\pi_1(\mathbb{A}^{1, \text{dR}})$
is via $\mathbb{A}^{1, \text{dR}} \to \pi_0(\mathbb{A}^{1, \text{dR}}) \xrightarrow{\mathrm{Frob}^i} \pi_0(\mathbb{A}^{1, \text{dR}})$.
To understand this group, let us utilize the cofiber sequence of
cotangent complexes from \cref{halftriangle} associated with
the diagram $W_n(k) \to k \to \mathbb{A}^{1, \text{dR}}$:
$$
\mathbb{L}_{k/W_n(k)} \otimes_k \mathbb{A}^{1, \text{dR}}
\to \mathbb{L}_{\mathbb{A}^{1, \text{dR}}/W_n(k)} \to \mathbb{L}_{\mathbb{A}^{1, \text{dR}}/k}.
$$
By \cref{crystalline}, the map $\w{End}_k (\mathbb{A}^{1, \w{dR}})\to \w{End}_k (\pi_0 (\mathbb{A}^{1, \w{dR}}))$ is injective with image $\underline{\w{Frob}}_k(B)$. Therefore, again by \cref{square-zero}, we have $\mathrm{Hom}_{\mathbb{A}^{1, \text{dR}}}(\mathbb{L}_{\mathbb{A}^{1, \text{dR}}/k}, \pi_1(\mathbb{A}^{1, \text{dR}})[1]) = 0$
and we get an injection
\[
\mathrm{Hom}_{\mathbb{A}^{1, \text{dR}}}(\mathbb{L}_{\mathbb{A}^{1, \text{dR}}/W_n(k)}, \pi_1(\mathbb{A}^{1, \text{dR}})[1])
\hookrightarrow 
\mathrm{Hom}_{\mathbb{A}^{1, \text{dR}}}(\mathbb{L}_{k/W_n(k)} \otimes_k \mathbb{A}^{1, \text{dR}}, \pi_1(\mathbb{A}^{1, \text{dR}})[1]),
\]
and the latter is identified with
\[
\mathrm{Hom}_k(\mathbb{L}_{k/W_n(k)}, \pi_1(\mathbb{A}^{1, \text{dR}})[1]) = \mathrm{Hom}_k(k[1], \pi_1(\mathbb{A}^{1, \text{dR}})[1])
= \pi_1(\mathbb{A}^{1, \text{dR}})(B)
= W[p](B).
\]
Here the first identification follows from the fact that $\tau_{\le 1}\mathbb{L}_{k/W_n(k)} = k[1]$,
and the last identification is due to \Cref{cohomology of A1dR}.(2).
Unraveling definitions, for any $u \in (1+ W[p])^\times(B)$,
the element $u \cdot \w{Frob}^i$ (see \Cref{construct aut} and \Cref{automorphism proposition}) in the fibre of $f_n$ is sent to $u-1 \in W[p](B)$.
One easily sees that the previous sentence in fact gives a bijection;
therefore, the fibre of $f_n$ over $\w{Frob}^i$ is exactly $(1+ W[p])^\times(B) \cdot \w{Frob}^i$ and finishes the calculation
of endomorphism sets.

The multiplication law is checked by chasing through the diagram:
on the quasi-ideal model, the homomorphism $u \cdot \w{Frob}^i$ sends an element $x \in W(B)$ to $u \cdot \w{Frob}^i(x)$,
and one computes $u \cdot \w{Frob}^i(v \cdot F^j(x)) = u \cdot \w{Frob}^i(v) \cdot \w{Frob}^{i+j}(x)$.
\end{proof}

\begin{remark}
In the above proof, one does not actually need to work with fpqc sheaves and the same proof works merely at the level of presheaves. However, if one only wanted to prove \cref{representing endomorphism monoid} in the case when $m=1,$ one can work with fpqc sheaves or quasi-syntomic sheaves and use the fact that $\pi_0 (\mathbb{A}^{1, \w{dR}})= \mathbb{G}_a$ from \cref{cohomology of A1dR} to simplify the proof and avoid invoking \cref{crystalline} and \cref{image is Frobenius powers}. The case $m=1$ is sufficient for our application in \cref{section5}.
\end{remark}{}

\begin{corollary}\label{corend}Let $k$ be an arbitrary perfect ring.
We consider the functor $\mathrm{End}_{m,n}$ from \cref{sec4} for a fixed $m \ge 1.$ There are natural maps of sheaves $\mathrm{End}_{m,n'} \to \mathrm{End}_{m,n}$ for $n' \ge n,$ which induces an isomorphism if $n \ge 2.$ If $n'>n,$ and $n=1,$ then all fibres of this natural map are given by the group scheme $(1+ W[p])^\times.$ The sheaf $\mathrm{End}_{m,1}$ is $\underline{\mathrm{Frob}}_k.$ 
\end{corollary}{}

\begin{proof}This follows from combining \cref{crystalline} and \cref{representing endomorphism monoid}.
\end{proof}{}
\begin{remark}
\leavevmode
\label{remark about End(dR)}
\begin{enumerate}

\item The stabilization of $\w{End}_{m,n}$ for $n \ge 2$ that we see above suggests that lifting to $W_n$ for $n>2$
poses no extra information on the de Rham cohomology of the special fibre, at least in a functorial sense.
In next section we will see the extra information on liftability to second Witt vectors gives
a strengthening to Deligne--Illusie's decomposition theorem \cite{DI87}.
Combining these two results,
we are led to believe the following dichotomy of possibilities on a follow-up question \cite[Remarques 2.6.(iii)]{DI87}:
either liftable over $W_2$ always guarantees Hodge--de Rham spectral sequence degenerates;
or there is a counterexample (necessarily of dimension $\geq p+1$) which is liftable all the way over $W$. \vspace{1mm}

\item 
If $B$ has characteristic $p$, then $p = V \circ F$ on $W(B)$.
The defining equation $u \cdot p = p$ of $(1+ W[p])^\times$ becomes
$V(F(u)) = V(1)$.
Since $V$ is always injective, the group scheme $(1+ W[p])^\times$ over a characteristic $p$ base
becomes $\mathbb{G}_m^{\sharp} \coloneqq W^{\times}[F]$,
namely the Frobenius kernel of the multiplicative group scheme $W^{\times}$. \vspace{1mm}

\item
The above discussion tells us that the functorial automorphism group scheme of the
mod $p$ de Rham cohomology theory on $W_2(k)$-algebras is given by $\mathbb{G}_m^{\sharp}$.
Note that there is a natural inclusion $\mu_p \to \mathbb{G}_m^\sharp$ which induces a product decomposition $\mathbb{G}_m^\sharp = \mu_p \times \mathbb{G}_a^\sharp$ (see \cref{productgm}). In \Cref{computing mup action} we will utilize the automorphisms coming from $\mu_p.$
The remaining $\mathbb{G}_a^{\sharp}$ worth of automorphisms is related to the Sen operator studied in \cite{BhaLur}.
\end{enumerate}
\end{remark}

Our calculation shows that there is no functorial splitting of the whole 
mod $p$ derived de Rham complex, as a functor from $W_2(k)$-algebras to $\mathrm{CAlg}(D(k))$,
into direct sums of the graded pieces of its conjugate filtrations.

\begin{proposition}
\label{no further functorial splitting}
There is no functorial splitting
\[
\mathrm{dR}_{(- \otimes_{W_2(k)} k)/k} \simeq
\bigoplus_{i \in \mathbb{N}_{\geq 0}} \mathrm{Gr}^{\mathrm{conj}}_{i}\big(\mathrm{dR}_{(- \otimes_{W_2(k)} k)/k}\big)
\]
as a functor from smooth $W_2(k)$-algebras to $\mathrm{CAlg}(D(k))$.
\end{proposition}

\begin{proof}
Indeed, if there were such a splitting, we would get an automorphism parametrized by $\mathbb{G}_m$,
with the $i$-th graded piece having pure weight $i$.
From the calculation of the endomorphism monoid in \cref{representing endomorphism monoid}, this would give us an injection 
$\mathbb{G}_m \hookrightarrow \mathbb{G}_m^{\sharp}$.
But the Frobenius on $\mathbb{G}_m$ is non-zero, whereas it is zero
on $\mathbb{G}_m^{\sharp}$, hence we know that there is no injective map
$\mathbb{G}_m \hookrightarrow \mathbb{G}_m^{\sharp}$ over any characteristic $p$ base,
therefore getting a contradiction.
\end{proof}

\begin{remark}[Twisted forms of de Rham cohomology]\label{twists}
\cref{representing endomorphism monoid} can be applied to understand a question considered by Antieau and Moulinos on possible existence of \'etale twists of the de Rham cohomology functor in some cases: let $k$ be a perfect ring and let $B$ be an ordinary $W_m(k)$-algebra, does there exist a functor $F: \w{ARings}_{W_n(k)} \to \w{CAlg}(D(B))$ which is isomorphic to $\w{dR}_{m,n} \otimes_{W_m(k)} B$ \'etale locally on $\w{Spec}\,B?$ We thank Antieau for mentioning this question to us. By \cref{representing endomorphism monoid}, such functors are classified by $\w{H}^1_{\textit{\'et}} (\w{Spec}\, B, (1+ W[p])^\times).$ When $m=1$ and $B$ is perfect, one can show that $\w{H}^1_{fpqc} (\w{Spec}\, B, (1+ W[p])^\times)=0$ by using $(1+ W[p])^\times \simeq \mathbb{G}_m^\sharp \simeq \mu_p \times \mathbb{G}_a^\sharp$ over $\w{Spec}\, B.$ So in that case, there does not even exist a non-trivial fpqc twist. However, the cohomology group can be non-zero for some choices of $B.$ It would be interesting to study the corresponding twisted forms of de Rham cohomology which can be seen as new cohomology theories, but that direction is not pursued further in this paper. It would also be interesting to compute $\w{H}^1_{\textit{\'et}}(\w{Spec}\, B, (1+ W[p])^\times)$ in general for $m>1.$
\end{remark}{}

\section{Application to the Deligne--Illusie decomposition}\label{section5}

\subsection{Drinfeld's refinement of the Deligne--Illusie decomposition}

In this section we will explain how to apply our result from \cref{representing endomorphism monoid} on endomorphisms of de Rham cohomology functor to recover a recent result of Drinfeld concerning a classical theorem
due to Deligne--Illusie \cite{DI87} and Achinger--Suh \cite{AchSuh}.

\begin{notation}
Fix a perfect ring $k$ as before and consider the monoid scheme $\w{End}_{1,n}$ from \cref{corend} over $k$. Let $B$ be a $k$-algebra and let $\sigma \in \mathrm{End}_{1,n}(B)$.
By definition we get an endomorphism induced by $\sigma$
\[
\mathrm{dR}_{R/W_n(k)} \otimes_{W_n(k)} B \xrightarrow{\sigma} \mathrm{dR}_{R/W_n(k)} \otimes_{W_n(k)} B,
\]
which is functorial in the $W_n(k)$-algebra $R$.
\end{notation}

\begin{definition}
For any $W_n(k)$-algebra $R$, we define the conjugate filtration $\Fil^{\conj}_i$ on
$\mathrm{dR}_{R/W_n(k)} \otimes_{W_n(k)} k$ to be the
left Kan extension of the canonical filtration on polynomial (or smooth) $W_n(k)$-algebras.
\end{definition}

\begin{lemma}
Assume that $k \to B$ is flat, then $\sigma$ preserves $\Fil^{\conj}_i \otimes_k B$ for all $i$.
\end{lemma}

\begin{proof}
By functoriality, any morphism must preserve the canonical filtration.
If $R$ is a polynomial (or smooth) $W_n(k)$-algebra,
one easily shows that the canonical filtration on $\mathrm{dR}_{R/W_n(k)} \otimes_{W_n(k)} B$
is just $\Fil^{\conj}_i \otimes_k B$.
\end{proof}


By \Cref{representing endomorphism monoid} and \Cref{remark about End(dR)} (4), we have an inclusion of
$k$-schemes $(\mathbb{G}_m^{\sharp}) \subset \mathrm{End}_{1,2}$.
Let $B = \mathrm{\Gamma}(\mathbb{G}_m^{\sharp}, \mathcal{O})$, then the identity map defines an element
$\sigma \in \mathbb{G}_m^{\sharp}(B)$, which can be regarded as the universal point.
By the above discussion, the universal point $\sigma$ gives rise to a comodule structure
on $\mathrm{dR}_{R/W_2(k)} \otimes_{W_2(k)} k$ over the Hopf algebra $B$, functorial in the $W_2(k)$-algebra
$R$ and the conjugate filtration is an increasing filtration of sub-comodules.
Alternatively, we may view this as an action of $\mathbb{G}_m^{\sharp}$ on the mod $p$ de Rham cohomology
$\mathrm{dR}_{-/W_2(k)} \otimes_{W_2(k)} k$.
We may ask the effect of $\mathbb{G}_m^{\sharp}$-action on each graded piece of the conjugate filtration,
viewed as a functor from the category of $W_2(k)$-algebras to the derived $\infty$-category
of $B$-comodules. The latter can be defined as the derived $\infty$-category of quasi-coherent sheaves on $B \mathbb{G}_m^\sharp.$


Recall that the category of $\mu_p$-representations is semi-simple with simple objects
given by $\mathbb{Z}/p$-worth of powers of the universal character.
We follow the convention that the universal character $\mu_p \hookrightarrow \mathbb{G}_m$
has weight $1$.
The following result is first observed by Drinfeld via prismatization, communicated to us by Bhatt.

\begin{theorem}
\label{computing mup action}
The action of $\mathbb{G}_m^{\sharp}$ on the $i$-th graded piece of the conjugate filtration
factors through the natural projection $\mathbb{G}_m^{\sharp} \to \mu_p$,
and the resulting $\mu_p$-action is of pure weight $i \in \mathbb{Z}/p$.
\end{theorem}
This fact also appears in \cite[Example 4.7.17]{BhaLur}, where it is
proved by contemplating the Sen operators. Below we give a different argument.
\begin{proof}
Derived Cartier isomorphism \cite[Proposition 3.5]{Bha12}
reduces us to showing the statement for $i = 0 \text{ and } 1$.
Since the conjugate filtration is defined via left Kan extension from its values
on polynomial algebras, using the classical Cartier isomorphism and K\"{u}nneth formula, 
we are reduced to understanding the behavior of $\sigma$ on the cohomology of
\[
\mathrm{dR}_{W_2(k)[x]/W_2(k)} \otimes_{W_2(k)} k 
\simeq \mathrm{dR}_{k[x]/k}.
\]
By observing that the whole situation is base changed from $k = \mathbb{F}_p$,
we immediately reduce to $k = \mathbb{F}_p$.

According to \Cref{unwinding construction}, the action of $\sigma$ is defined via
the identification
\[
\mathrm{dR}_{\mathbb{Z}/p^2[x]/(\mathbb{Z}/p^2)} \otimes_{\mathbb{Z}/p^2} B
\simeq \mathrm{R\Gamma}(\mathbb{A}^{1,\mathrm{dR}}_{B}, \mathcal{O})
\]
and the homomorphism of $\mathbb{Z}/p^2$-algebra stack over $\mathbb{G}_m^{\sharp}$ given by the diagram below
(here $W_{B}$ denotes the Witt ring scheme over the base scheme $\mathbb{G}_m^{\sharp}$)
\[
\xymatrix{
W_{B} \ar[r]^{\times p} \ar[d]_{\times \sigma} & W_{B} \ar[d]^{\text{id}} \\
W_{B} \ar[r]^{\times p} & W_{B}.
}
\]

The first cohomology of
$\mathrm{dR}_{{\mathbb{F}_p}[x]/{\mathbb{F}_p}}$ is a free rank $1$ module over its zeroth cohomology.
Therefore, all we need to do is the following.
\begin{enumerate}
\item Show that the induced map on
$\w{H}^0(\mathbb{A}^{1,\mathrm{dR}}_{B}, \mathcal{O})$ is trivial.
\item Exhibit a non-zero element 
$v \in \w{H}^1(\mathbb{A}^{1,\mathrm{dR}}_{\mathbb{F}_p}, \mathcal{O})$
which pulls back to a weight $1$ element in $\w{H}^1(\mathbb{A}^{1,\mathrm{dR}}_{B}, \mathcal{O})$.
\end{enumerate}
To avoid confusion, in the quasi-ideal $W \xrightarrow[]{\times p} W,$ let us denote the ring scheme on the right as $W \coloneqq \mathrm{Spec}({\mathbb{F}_p}[X_0, X_1, \ldots])$
and the $W$ on the left as $W\coloneqq\mathrm{Spec}({\mathbb{F}_p}[Y_0, Y_1, \ldots])$.
Here $X_i$'s (and similarly $Y_i$'s) are the Witt coordinates.
In the above diagram, one easily checks the effect of $\text{id}^*$ and $(\times \sigma)^*$ on the following elements:
$X_i \mapsto X_i$, $Y_0 \mapsto t_0 \cdot Y_0$.
Here $t_0$ denotes the element in $B$ corresponding to the natural projection $\mathbb{G}_m^{\sharp} \to \mu_p$.


Now (1) is easily verified: $\w{H}^0(\mathbb{A}^{1,\mathrm{dR}}_{B}, \mathcal{O}) \subset B[X_0, X_1, \ldots]$,
hence invariant under the $\sigma$-action.

As for (2), we claim that $1 \otimes Y_0 \in \mathbb{F}_p[X_i, Y_j]$ 
is a nonzero class in $\w{H}^1(\mathbb{A}^{1,\mathrm{dR}}_{\mathbb{F}_p}, \mathcal{O})$.
Here we are using the \v{C}ech nerve of $\mathrm{Spec}({\mathbb{F}_p}[X_0, X_1, \ldots]) \to \mathbb{A}^{1,\mathrm{dR}}_{\mathbb{F}_p}$
to calculate the cohomology of $\mathbb{A}^{1,\mathrm{dR}}_{\mathbb{F}_p}$,
implicitly we have used the fact that the $[1]$-term of the \v{C}ech nerve is given by $\mathrm{Spec}({\mathbb{F}_p}[X_i, Y_j])$.
Granting this claim, the action of $\sigma$ sends $Y_0$ to $t_0 \cdot Y_0$, hence the action on the
class $1 \otimes Y_0$ is via the natural projection $\mathbb{G}_m^{\sharp} \to \mu_p$ and has weight $1$.
To prove the claim, we use the maps
\[
W[F] = \mathbb{G}_a^{\sharp} \to W = \mathrm{Spec}({\mathbb{F}_p}[Y_0, Y_1, \ldots]) \to \mathbb{G}_a = \mathrm{Spec}({\mathbb{F}_p}[Y_0]),
\]
where the middle $W$ is the copy of quasi-ideal $W$.
The above maps induce a sequence of abelian group stacks:
\[
B \mathbb{G}_a^{\sharp} \to \mathbb{A}^{1,\mathrm{dR}} \to B \mathbb{G}_a.
\]
Recall that there is a canonical identification $\w{H}^1(BG, \mathcal{O}) \simeq \mathrm{Hom}(G, \mathbb{G}_a)$ for affine group schemes $G$ via faithfully flat descent along $* \to BG.$
The identity map on $\mathbb{G}_a = \mathrm{Spec}({\mathbb{F}_p}[Y_0])$ pulls back to $1 \otimes Y_0 \in \mathbb{F}_p[X_i, X_j]$,
which checks that $1 \otimes Y_0$ is a cocycle.
Furthermore, recall that the induced map $\mathbb{G}_a^{\sharp} \to \mathbb{G}_a$ realizes the former as the
divided power envelope of the origin inside the latter, in particular it is a nonzero map.
From the above identification, this tells us that $1 \otimes Y_0$ pulls to a nonzero class in
$\w{H}^1(B\mathbb{G}_a^{\sharp}, \mathcal{O})$.
Therefore, the class $1 \otimes Y_0$ is a nonzero class in $\w{H}^1(\mathbb{A}^{1,\mathrm{dR}}_{\mathbb{F}_p}, \mathcal{O})$.
\end{proof}

\begin{remark}
Let us mention another way to obtain the above result concerning the
$\mathbb{G}_m^{\sharp}$-action on $\mathrm{dR}_{k[x]/k}$.
As explained, the action arises from the action of $\mathbb{G}_m^{\sharp}$
on $\mathbb{A}^{1,\mathrm{dR}}$ in characteristic $p$. One can show that the stack underlying $\mathbb{A}^{1,\w{dR}}_{\mathbb{F}_p}$ (without the ring stack structure) decomposes as $\mathbb{G}_a \times B\mathbb{G}_a^{\sharp}$ (\textit{c.f.}~\cite[Proposition 5.12 ]{BhaLur2}) and the action of $\mathbb{G}_m^{\sharp}$ is trivial on $\mathbb{G}_a$
and weight $1$ on $B\mathbb{G}_a^{\sharp}.$ This gives the desired statement.
We thank the referee for pointing this out to us.
\end{remark}

Note that the natural projection $\mathbb{G}_m^{\sharp} \to \mu_p$ admits a splitting:
the Teichm\"{u}ller lift defines a map of group schemes
$\mathbb{G}_m \hookrightarrow W^{\times}$, which induces a map of group schemes
$\mu_p \hookrightarrow \mathbb{G}_m^{\sharp}$.

Let $X$ be a smooth scheme over $W_2(k)$, consider the de Rham cohomology of its special fibre (relative to $k$)
which by the above discussion admits a $\mu_p$-action.
Now, if we look at the canonical truncation in a range of width at most $p$,
then the weights that show up in $\mathbb{Z}/p$ are pairwise distinct, hence we get a splitting
of the induced conjugate filtration.
Therefore, the above theorem implies the following improvement of a result due to Achinger and Suh \cite[Theorem 1.1]{AchSuh}, which in turn
is a strengthening of Deligne--Illusie's result \cite[Corollaire 2.4]{DI87}. 

\begin{corollary}[Drinfeld]
\label{Drinfeld splitting}
Let $k$ be a perfect ring of characteristic $p > 0$, let $X$ be a smooth scheme over $W_2(k)$,
and let $a \leq b \leq a + p-1$.
Then the canonical truncation $\tau_{[a,b]}(\Omega^{\bullet}_{X_k/k})$ splits.
\end{corollary}

Note that when $p>2$, in Achinger--Suh's statement in loc.~cit.~they need $b < a + p - 1$ so their allowed width needs to be at most $p-1$.
In fact, more generally, we have the following decomposition as a consequence of the $\mathbb{G}_m^{\sharp}$-action in described in \Cref{computing mup action}.

\begin{corollary}[Drinfeld, {\textit{c.f.}~\cite[Remark A.5]{AchSuh}}]
\leavevmode
Let $X$ be a smooth scheme over $W_2(k)$ with special fibre $X_k$.
Then there exists a splitting, functorial in $X$, in the derived $\infty$-category
of Zariski sheaves on $X'_k$ 
\[
F_{X_k/k,*}(\mathrm{dR}_{X_k/k}) \simeq \oplus_{i \in \mathbb{Z}/p} F_{X_k/k,*}(\mathrm{dR}_{X_k/k}^{\mathrm{weight }=i}).
\]
Moreover $\mathcal{H}^j(F_{X_k/k,*}(\mathrm{dR}_{X_k/k}^{\mathrm{weight }=i})) \not= 0$ implies $j \equiv i$ in $\mathbb{Z}/p$.
Here $X'_k$ is the Frobenius twist of $X_k$ and $F_{X_k/k}$ is the relative Frobenius.
In particular, the conjugate spectral sequence of liftable smooth varieties can have
non-zero differentials only on $(mp+1)$-st pages, where $m \in \mathbb{Z}_{> 0}$.
\end{corollary}

\begin{remark}
Drinfeld observed the results in this subsection by using the ``stacky approach'' to prismatic crystals 
(which he calls ``prismatization''),
which is independently developed by Bhatt--Lurie \cite{BhaLur} as well. 
Using the prismatization functor, Drinfeld produced an action of $\mu_p$ on the de Rham complex of a smooth scheme over $k$ that lifts to $W_2(k)$.
Our paper partly grew out of an attempt of making sense of and reproving Drinfeld's theorem
without introducing prismatization and taking a very algebraic/categorical approach instead. 
In the recent paper of Bhatt--Lurie \cite{BhaLur}, 
this action is obtained in a more geometric way by understanding the prismatization of $\w{Spec}(W_2(k)).$
\end{remark}

\subsection{Uniqueness of functorial splittings} \Cref{Drinfeld splitting} provides a functorial splitting of the $(p-1)$-st conjugate
filtration of the mod $p$ derived de Rham cohomology of any $W_2(k)$-algebra.
On the other hand, the classical Deligne--Illusie splitting also has an $\infty$-categorical functorial
enhancement~\cite[Theorem 1.3.21 and Proposition 1.3.22]{KP19} which, in spirit, is
more related to the work of Fontaine--Messing \cite{FM87} and Kato \cite{Ka87}.

It is a natural question to ask whether these two splittings agree in a functorial way.
By how these two splittings are defined, we see immediately that they are both compatible with
the module structure over the $0$-th conjugate filtration and induced from the splitting of the
$1$-st conjugate filtration by an averaging process, \textit{c.f.}~the step (a) in proof of \cite[Th\'{e}or\`{e}me 2.1]{DI87}.

Below we will prove that there is a unique way to functorially split the first conjugate filtration,
hence the above two functorial splittings must be the same.
To that end, let us fix some notations.
\begin{notation}
Let us consider the stable $\infty$-category $\w{Fun}(\w{Alg}_{W_2(k)}^{\w{sm}}, D(k)),$ where $\w{Alg}_{W_2(k)}^{\w{sm}}$ is the category of smooth $W_2(k)$-algebras and $D(k)$ is the derived (stable) $\infty$-category of $k$-vector spaces. Let us use $\mathcal{O}$ to denote the functor that sends any $W_2(k)$-algebra $R$
to the $0$-th conjugate filtration of $\mathrm{dR}_{(R \otimes_{W_2(k)} k)/k}$, which has the structure of a commutative algebra object in $\w{Fun}(\w{Alg}_{W_2(k)}^{\w{sm}}, D(k)).$ The functor obtained by considering the $1$-st piece of the conjugate filtration will be denoted $M$,
viewed as an $\mathcal{O}$-module.
We have a natural map $\mathcal{O} \to M$; we denote the cofiber by $G$, which is the first graded piece of the conjugate filtration, also viewed as an $\mathcal{O}$-module.
\end{notation}

Now we have a cofiber sequence of $\mathcal{O}$-modules: $\mathcal{O} \to M \to G$.

\begin{theorem}
\label{uniqueness of functorial splitting of Conj1}
In the above notations, there is a unique functorial $\mathcal{O}$-module splitting
\[
M = \mathcal{O} \oplus G
\] in $\w{Fun}(\mathrm{Alg}_{W_2(k)}^{\w{sm}}, D(k)).$
In particular the splitting of $\Fil^{\mathrm{conj}}_{p-1}(\mathrm{dR}_{(- \otimes_{W_2(k)} k)/k})$
obtained in \Cref{Drinfeld splitting} and \cite[Theorem 1.3.21]{KP19} agree.
\end{theorem}

\begin{proof}
The existence part is provided by either \Cref{Drinfeld splitting} or \cite[Theorem 1.3.21]{KP19}.
We focus on the uniqueness part in this proof.

Firstly, we note that it suffices to show the uniqueness of splitting
as a quasi-syntomic sheaf on the quasi-syntomic site of $W_2(k)$.
This is because they are left Kan extended from the polynomial case,
and polynomial algebras are quasi-syntomic.
The site $\qsyn_{W_2(k)}$ admits a basis of large quasi-syntomic $W_2(k)$-algebras, so we may restrict our functors to this subclass of $W_2(k)$-algebras and show the uniqueness 
of splitting there. 
All three functors have discrete values on this subclass of $W_2(k)$-algebras,
so $\mathcal{O}$ is a sheaf of ordinary $k$-algebras given by
$R \mapsto R/p$ (up to Frobenius twist) and $M$ and $G$ are sheaves of ordinary $\mathcal{O}$-modules.
We will show that there exists a unique section to the surjection of sheaves
of $\mathcal{O}$-modules $M \twoheadrightarrow G$.

Step 1: Let us consider the algebra $R = W_2(k)[x^{1/p^{\infty}}]/(x)$.
In this case, we have that $$D \coloneqq \mathrm{dR}_{(R \otimes_{W_2(k)} k)/k} \simeq D_{(x)}(k[x^{1/p^{\infty}}])$$
is the divided power envelope of $(x)$ in $k[x^{1/p^{\infty}}]$.
This algebra admits a natural grading by the monoid $\mathbb{N}[1/p]$.
The values of our sheaves evaluated at $R$ are $\mathcal{O} = k[x^{1/p^{\infty}}]/(x^p)$,
and $M$ is the degree $[0, 2p)$ part of $D_{(x)}(k[x^{1/p^{\infty}}])$,
whereas $G$ is the degree $[p, 2p)$ part.
One checks easily that $G$ is generated by $\overline{\gamma_p(x)}$ (mod degree $[0,p)$ part)
as an $\mathcal{O}$-module in this case.
We claim that the section necessarily sends this generator to $\gamma_p(x) \in M$.
Suppose that the section sends this generator to some element $f(x) \in M$.
We look at the following two maps of $W_2(k)$-algebras $R \to R \otimes_{W_2(k)} W_2(k)[t^{1/p^{\infty}}]$,
with $x^m \mapsto x^m \cdot t^m$ and $x^m \mapsto x^m$.
The associated mod $p$ derived de Rham cohomology is given by $D \otimes_k k[t^{1/p^{\infty}}]$.
Since the corresponding map of values of $G$ is
\[
\overline{\gamma_p(x)} \mapsto \overline{\gamma_p(tx)} = t^p \overline{\gamma_p(x)}
\text{ and } 
\overline{\gamma_p(x)} \mapsto \overline{\gamma_p(x)},
\]
functoriality tells us
that $t^p f(x) = f(tx) \in D \otimes_k k[t^{1/p^{\infty}}]$.
This implies that $f(x)$ is a homogeneous degree $p$ element in $M$ which maps to $\overline{\gamma_p(x)} \in G$;
therefore, it must be $\gamma_p(x) \in M$.

Step 2: Next we consider the algebra $R_n = W_2(k)[x_i^{1/p^{\infty}}; i = 1, \ldots, n]/(\sum_{i=1}^n x_i)$.
In this case, we have
\[
D_n \coloneqq \mathrm{dR}_{(R_n \otimes_{W_2(k)} k)/k} = D_{(\sum_{i=1}^n x_i)}(k[x_i^{1/p^{\infty}}]).
\]
Then the values of our sheaves evaluated at $R_n$ is given by $\mathcal{O} = k[x_i^{1/p^{\infty}}]/(\sum_{i=1}^n x_i^p)$,
and $M = \mathcal{O} \cdot \{1, \gamma_p(\sum_{i=1}^n x_i)\}$
whereas $G = \mathcal{O} \cdot \overline{\gamma_p(\sum_{i=1}^n x_i)}$.
In this case, we claim that the section necessarily sends $\overline{\gamma_p(\sum_{i=1}^n x_i)}$
to $\sum_{i=1}^n \gamma_p(x_i)$.
Note that this sum makes sense as an element in $D_n$ and in fact is in $M$,
for instance, one may repeatedly use $\gamma_p(x+y) = \sum_{i= 0}^p \gamma_i(x) \cdot \gamma_{p-i}(y)$
to see this.
Now, to see the above claim, we first use the same argument as in the previous paragraph
to see that the section of $\overline{\gamma_p(\sum_{i=1}^n x_i)}$ is necessarily
a homogeneous degree $p$ element $f(x_i)$.
Then we use the functoriality provided by the map
$R_n \to R^{\otimes_{W_2(k)} n} = W_2(k)[x_i^{1/p^{\infty}}; i = 1, \ldots, n]/(x_i; i = 1, \ldots, n)$
to see that the element $g(x_i) \coloneqq f(x_i) - \sum_{i=1}^n \gamma_p(x_i)$
is a homogeneous degree $p$ element in the kernel of the induced map
$D_n \to D^{\otimes_k n}$.
The degree $p$ part of the kernel is the $k$-span of
$\{x_i^p\}_{i = 1}^n$ modulo $k \cdot \sum_{i=1}^n x_i^p$.
Finally, using functoriality with respect to switching variables,
we see that $g(x_i)$ must be a permutation-invariant element,
hence necessarily $0$ unless $n = p = 2$.
Therefore, we know that when $n \geq 3$, the associated section is determined uniquely.
By functoriality, the section associated with $R_3$ determines the section associated with $R_2$,
this finishes the proof of our claim above.

Step 3: The universal algebra that we need to consider now is 
$R' = W_2(k)[x^{1/p^{\infty}}, y^{1/p^{\infty}}]/(x + py)$.
Note that $R'/p = R/p \otimes_k k[y^{1/p^{\infty}}]$, so the value of relevant sheaves
are those in Step 1 tensored over $k$ with $k[y^{1/p^{\infty}}]$.
The generator $\overline{\gamma_p(x)} = \overline{\gamma_p(x+py)}$ of $G$
under a functorial section goes to $\gamma_p(x) + g(x,y)$,
where $g(x,y) \in k[x^{1/p^{\infty}}, y^{1/p^{\infty}}]/(x^p)$ has degree $p$
by the same argument as in Step 1.
We claim $g(x,y) = \frac{y^p}{(p-1)!}$.
To see this, we first observe that
\[
x_1 + x_2 = (x_1^{1/p} + x_2^{1/p})^p + p \cdot F(x_1, x_2) \text{ in }W_2(k)[x_1^{1/p}, x_2^{1/p}],
\]
where we view $F(x_1, x_2) \in k[x^{1/p^{\infty}}, y^{1/p^{\infty}}]$, a degree $1$ polynomial.
Then we see that there is a map $R' \to R_2$ sending $x$ and $y$ to Teichm\"{u}ller lifts of
$x_1 + x_2$ and $F(x_1, x_2)$.
The induced map of corresponding $D$'s sends $\gamma_p(x) + g(x,y)$ to
$\gamma_p(x_1 + x_2) + g(x_1 + x_2, F(x_1, x_2))$.
On the other hand the functoriality forces this element to be sent to
$\gamma_p(x_1) + \gamma_p(x_2)$ by Step 2.
Therefore, we get a relation
\[
\gamma_p(x_1 + x_2) + g(x_1 + x_2, F(x_1, x_2)) = \gamma_p(x_1) + \gamma_p(x_2).
\]
Let $h(x,y) = g(x,y) - \frac{y^p}{(p-1)!} \in k[x^{1/p^{\infty}}, y^{1/p^{\infty}}]/(x^p)$
which also has degree $p$.
Combining these relations, we see that $h(x_1 + x_2, F(x_1, x_2)) = 0 \in k[x_1^{1/p^{\infty}}, x_2^{1/p^{\infty}}]/(x_1^p + x_2^p)$.
Applying the next lemma with $x_1 + x_2 = a$ and $x_2 = b$, we conclude that $h(x,y)$ must be $0$.

Step 4: We finish the proof in this step.
Given any large quasi-syntomic $W_2(k)$-algebra $S$, we can find an algebra
$S'$ of the form $W_2(k)[X_i^{1/p^\infty}, Y_j^{1/p^\infty}; i \in I, j \in J]/(Y_j + f_j(X_i); j \in J)$
and a surjection $S' \twoheadrightarrow S$ inducing a surjection of their values
of all the relevant sheaves, see the proof of \cite[Proposition 7.10]{BS19} or \cite[Theorem 3.13]{LL20}
for details.
The value of $G$ in this case is generated, as an $\mathcal{O}$ module, by $\gamma_p(Y_j + f_j(X_i))$
where $j \in J$.
By functoriality, we may reduce to the case where 
$S = W_2(k)[\underline{X}^{1/p^\infty}, Y^{1/p^\infty}]/(Y + f(\underline{X}))$.
In this case $G$ is generated over $\mathcal{O}$ by the element $\overline{\gamma_p(Y + f(\underline{X}))}$
we want to show that the section is forced on this element. Note that any element in $W_2(k)[\underline{X}^{1/p^\infty}, Y^{1/p^\infty}]$ can be written as
$[P_1] + p \cdot [P_2]$, a Teichm\"{u}ller lift plus $p$ times another Teichm\"{u}ller lift.
Therefore, we can define a map $R' \to S$ sending $X$ to $[P_1]$ and $Y$ to $[P_2]$.
Then we see that the section of $\overline{\gamma_p(Y + f(\underline{X}))}$
must be $\gamma_p(P_1) + \frac{P_2^p}{(p-1)!}$ by Step 3.
This shows the rigidity as desired.
\end{proof}

\begin{lemma}
Let $F(a,b) \in k[a^{1/p^\infty}, b^{1/p^\infty}]$ be the degree $1$ element such that its lift $\widetilde{F}$ to 
$W_2(k)[a^{1/p^\infty}, b^{1/p^\infty}]$ satisfies
\[
(a-b)^p + b^p = a^p + p \cdot \widetilde{F}(a^p, b^p) \text{ in } W_2(k)[a^{1/p^\infty}, b^{1/p^\infty}].
\]
Let $H(a, b) \in k[a^{1/p^\infty}, b^{1/p^\infty}]$ be a degree $p$ element which does not contain
the term $a^p$.
Suppose $H(a, F(a,b)) \in k[a^{1/p^\infty}, b^{1/p^\infty}]$ is divisible by $a^p$, then $H(a,b) = 0$.
\end{lemma}

\begin{proof}
Note that $F(a,b) = \sum_{i = 1}^{p-1} c_i \cdot a^{i/p} b^{(p-i)/p}$ with all of $c_i \not = 0$. The $a$-degree of $F(a,b)$ is less than $1$; therefore, the $a$-degree of $H(a, F(a,b))$ must be smaller than $p$
unless $H(a, F(a,b)) = 0$
(as $H(a,b)$ does not contain $a^p$ term).
The $a^p$ divisibility now forces $H(a, F(a,b)) = 0$.
Contemplating the $b$-degree of $H(a, F(a,b))$ shows that in fact $H(a,b)$ has to be $0$ to begin with.
\end{proof}

In the step 3 above, one can alternatively argue using the map $R_{p+1} \to R'$
sending $x_1$ to $x$ and the rest of the $p$ variables to $y$.


\appendix
\section{Topos theoretic cotangent complex}
\label{appendix}
The theory of cotangent complex appears in many places in the literature. For example, it has been discussed in \cite{Ill} in the context of simplicial ring objects in a $1$-topos and in \cite{Lur2}, where an $\infty$-categorical theory has been discussed for animated ring objects in spaces. However, in the proof of \cref{representing endomorphism monoid}, we required a formalism of cotangent complex in the generality of animated ring objects in an $\infty$-topos. In this appendix, we will sketch a formalism of cotangent complex in the above generality and its very basic properties which is sufficient for the proof of \cref{representing endomorphism monoid}. Our exposition basically uses the techniques from \cite{Lur2} and lifts them to the generality we need.

For simplicity, we will focus on the case necessary for our application where the $\infty$-topos $\mathcal{X}$ arises as sheaves of spaces on some Grothendieck site $\mathcal{C}$ which will be fixed. As in \cref{anim}, one defines the $\infty$-category $\w{ARings}(\mathcal{X}):=\w{ARings}(\mathcal{X})_{\mathbb{Z}}$ which is equivalent to the $\infty$-category of sheaves of animated rings on $\mathcal{C}.$ For a fixed animated ring $B$ in $\mathcal X$, one can also consider the $\infty$-category of connective $B$-modules in $\mathcal{X}$ defined as the category of sheaves on $\mathcal{C}$ (with values in animated abelian groups) of $B$-modules.

For $n \ge 0$, an object $F \in \w{ARings}(\mathcal{X})$ will be called $n$-truncated if $F(c)$ is $n$-truncated (i.e.,~$\pi_i(F(c)) = 0$ for all $i > n$)
for all $c \in \mathcal{C}.$ We let $\tau_{\le n} \w{ARings}(\mathcal{X}) \to \w{ARings}(\mathcal{X})$ denote the inclusion of the full subcategory of $n$-truncated objects in $\w{ARings}(\mathcal{X}).$ This admits a left adjoint that sends $G \mapsto \tau_{\le n} G$ which is obtained by $n$-truncating $G$ as a presheaf first and then applying sheafification.

\begin{construction}[The cotangent complex]\label{cotangent}Let $A \to B$ be a map in $\w{ARings}(\mathcal{X}).$ 
For any connective $B$-module $M,$ one can form the trivial square zero extension $B \oplus M,$ which is an object of $\w{ARings}(\mathcal{X})_A.$ There is a natural projection map $B \oplus M \to B,$ which regards $B \oplus M$ as an object of $(\w{ARings}(\mathcal{X})_A)_{/B}.$ One can consider the functor $M \to \w{Maps}_{(\w{ARings}(\mathcal{X})_A)_{/B}} (B, B \oplus M).$ By the adjoint functor theorem, this functor is corepresented by a connective $B$-module, which we will denote as $\mathbb{L}_{B/A}.$

\end{construction}{}

\begin{remark}Let $A \to B$ be a map in $\w{ARings}(\mathcal{X}).$ It follows that $\mathbb{L}_{B/A}$ defined as above is sheafification of the presheaf on $\mathcal{C}$ with values in animated abelian groups that sends $c \mapsto \mathbb{L}_{B(c)/A(c)}$ for $c \in \mathcal{C}.$ It naturally inherits the structure of a sheaf of connective $B$-modules on $\mathcal{C}.$\end{remark}

\begin{proposition}\label{halftriangle} For a sequence of morphisms 
$A \to B \to C$ in $\mathrm{ARings}(\mathcal{X})$, we have a cofiber sequence $$\mathbb{L}_{B/A} \otimes_B C \to \mathbb{L}_{C/A} \to \mathbb{L}_{C/B}$$ in the $\infty$-category of connective $C$-modules.
\end{proposition}{}
\begin{proof}This follows from \cref{cotangent}.
\end{proof}{}

\begin{remark}Let $C \in \mathrm{ARings}(\mathcal{X}).$ Let $U \to V \to W$ be a cofiber sequence in the $\infty$-category of connective $C$-modules. For any connective $C$-module $M,$ we obtain a long exact sequence 
$$\ldots \to \pi_1 \w{Maps} (W, M) \to \pi_1 \w{Maps} (V, M) \to \pi_1 \w{Maps} (U, M) \to \pi_0 \w{Maps} (W, M) \to \pi_0 \w{Maps} (V, M) \to \pi_0 \w{Maps} (U, M).$$
\end{remark}{}

\begin{definition}[Square-zero extensions] Let $A \in \w{ARings}(\mathcal X)$ and $B \in \w{ARings}(X)_A.$ Let $M$ be a connective $B$-module. A square-zero extension of $B$ by $M$ will be classified by $\w{Maps}_B(\mathbb{L}_{B/A}, M[1]),$ where the maps are considered in the $\infty$-category of connective $B$-modules. By \cref{cotangent}, square-zero extensions can be equivalently classified by $\w{Maps}_{(\w{ARings}(\mathcal{X})_A)_{/B}} (B, B \oplus M[1]).$ Given $s: B \to B\oplus M[1]$ which gives a section to the projection, the pullback $ B':=B \otimes_{B \oplus M[1]} B$ recovers the total space of the square-zero extension, where $B$ maps to $B \oplus M[1]$ via $s$ and the zero section. The fibre of $B' \to B$ can be identified with $M$ with the natural structure of an $A$-module.
\end{definition}{}

\begin{proposition}\label{square-zero}Let $C \in \mathrm{ARings}(\mathcal{X})_A$ and let $B' \to B$ in $\mathrm{ARings}(\mathcal{X})_A$ be a square-zero extension of $B$ by a connective $B$-module $M.$ 
There is a natural map $\mathrm{Maps}_{\mathrm{ARings}(\mathcal{X})_A}(C, B' ) \to \mathrm{Maps}_{\mathrm{ARings}(\mathcal{X})_A}(C, B)$ such that the non-empty fibres are torsors under the group $\mathrm{Maps}_C(\mathbb{L}_{C/A}, M),$ where the maps are taken in the category of connective $C$-modules. The $C$-module structure on $M$ is obtained via the map $C \to B$ over which the fiber is being taken.
\end{proposition}{}

\begin{proof}
Unwrapping the definitions and using the fact that the mapping spaces are $\infty$-groupoids, one can reduce to checking this in the case when $B'= B \oplus M$ is the trivial square zero extension. In this case, we can fix a map $C \to B$ and then we need to show that 
$\mathrm{Maps}_{(\mathrm{ARings}(\mathcal{X})_A)/B}(C, B \oplus M )$ is equivalent to $\w{Maps}_C  (\mathbb{L}_{C/A}, M).$ For this, we note that pullback along $C \to B$ gives an equivalence $\mathrm{Maps}_{(\mathrm{ARings}(\mathcal{X})_A)/B}(C, B \oplus M ) \simeq \mathrm{Maps}_{(\mathrm{ARings}(\mathcal{X})_A)/C}(C, C \oplus M ).$ By definition, $\mathrm{Maps}_{(\mathrm{ARings}(\mathcal{X})_A)/C}(C, C \oplus M ) \simeq \w{Maps}_C  (\mathbb{L}_{C/A}, M),$ which gives the conclusion.\end{proof}

\begin{remark}
For any object $A \in \mathrm{ARings}(\mathcal{X})_A,$ one can use the truncation functors to build a sequence of square-zero extensions $\ldots \tau_{\le n+1} A \to \tau_{\le n} A \to \tau_{\le n-1} A \ldots \to \tau_{\le 0} A= \pi_0(A).$ This can be seen by first showing a similar statement at the presheaf level and then sheafifying; at the presheaf level, the statement follows from the analogous statement for animated rings, \textit{c.f.}~\cite[Proposition 3.3.6]{Lur2}. In particular, if $A \in \mathrm{ARings}(\mathcal{X})_A $ is $1$-truncated, then $A$ is a square zero extension of $\pi_0 (A)$ by $\pi_1 (A) [1],$ where the latter is viewed as a connective $\pi_0(A)$-module in $\mathcal{X}.$
\end{remark}

\section{A product formula for $(1+ W[p])^\times$ in char.~$p>0$}
\label{productgm}
The group scheme $(1+ W[p])^\times$ was defined in \cref{notation of aut gp}. We will work over a fixed base ring of char.~$p>0$; then this group scheme is isomorphic to $W^\times [F],$ see \cref{remark about End(dR)} (2). The following proposition was stated in \cite[Lemma 3.3.4]{prismatization} and a more general proposition over $\mathbb{Z}_p$ has been proven in \cite[Proposition B.5.6]{drinew} by using the logarithm constructed in loc.~cit.;
see also \cite[Lemma 3.5.18]{BhaLur}. 
Let us give a more direct argument in char.~$p$ that does not use the logarithm and is closer to deformation theory in spirit.  
\begin{proposition}
There is a natural isomorphism $W^\times[F] \simeq W[F] \times \mu_p$ over any base ring of char.~$p.$
\end{proposition}{}

\begin{proof} Note that given any non-unital ring $(I, +, \cdot),$ one can define a monoid associated to it, which will be denoted as $I'.$ At the level of underlying sets, $I^{'} := I,$ but the composition $x * y$ is defined to be $x+ y + x \cdot y.$ Using the above construction along with the Yoneda lemma produces a functor from the category of non-unital ring schemes (e.g.,~ideals in unital ring schemes) to the category of monoid schemes. Note that we have a short exact sequence $$0 \to W[F] \to W[F] \xrightarrow[]{f} \alpha_p \to 0 $$ of group schemes. Moreover, the map $f: W[F] \to \alpha_p$ is a map of non-unital ring schemes when $W[F]$ and $\alpha_p \simeq \mathbb{G}_a[F]$ are both equipped with their natural non-unital ring scheme structures. Applying the functor we constructed before, we obtain a map $f': W^\times [F] \to \mu_p.$ It is clear that $f'$ is surjective.
The map $f'$ can be identified with projection to $0$-th Witt coordinate: given any test algebra $S$ and an element $x \in W[F](S)$,
then $f'$ sends $1 + x$ to $1 + x_0$ where $x_0$ is the $0$-th Witt coordinate.
In particular we see that the map $\mu_p \to W^{\times}[F]$ given by Teichm\"{u}ller lift is a section to $f'$.
It remains to identify $\w{Ker} f'$ with $W[F]$  as a group scheme. This follows from the lemma below.

\begin{lemma} For the map $f: W[F] \to \alpha_p,$ the ideal $\mathrm{Ker} f$ is a square-zero ideal.
\end{lemma}{}
\begin{proof}
We note that the multiplication in $W[F]$ is inherited from the ring scheme $W.$ Let $S$ be a test algebra of char.~$p$ and let $m, n \in (\w{Ker} f) (S).$ Then $m= V(m')$ and $n = V(n')$ for some $m', n' \in W(S).$ Here $V$ denotes the Verschiebung operator. We have $m\cdot n = V(m') \cdot n= V(m' \cdot F(n))=0,$ since $F(n) = 0.$
\end{proof}{}
The proposition now follows since we obtain a split exact sequence of group schemes $$0 \to W[F] \to W^{\times} [F] \xrightarrow[]{f'} \mu_p \to 0.$$\end{proof}{}

\bibliographystyle{amsalpha}
\bibliography{endobib}

\providecommand{\bysame}{\leavevmode\hbox to3em{\hrulefill}\thinspace}
\providecommand{\MR}{\relax\ifhmode\unskip\space\fi MR }
\providecommand{\MRhref}[2]{%
  \href{http://www.ams.org/mathscinet-getitem?mr=#1}{#2}
}
\providecommand{\href}[2]{#2}
\begin{thebibliography}{{Sta}21}

\bibitem[AS21]{AchSuh}
Piotr Achinger and Junecue Suh, \emph{Some refinements of the
  {Deligne}--{Illusie} theorem}, 2021, arXiv:2003.09857.

\bibitem[Bha12]{Bha12}
Bhargav Bhatt, \emph{p-adic derived de {Rham} cohomology}, 2012,
  arXiv:1204.6560.

\bibitem[BL22a]{BhaLur}
Bhargav Bhatt and Jacob Lurie, \emph{Absolute prismatic cohomology}, 2022,
  arXiv:2201.06120.

\bibitem[BL22b]{BhaLur2}
\bysame, \emph{The prismatization of $p$-adic formal schemes}, 2022,
  arXiv:2201.06124.

\bibitem[BLM21]{BLM20}
Bhargav Bhatt, Jacob Lurie, and Akhil Mathew, \emph{Revisiting the de
  {R}ham-{W}itt complex}, Ast\'{e}risque (2021), no.~424, viii+165.
  \MR{4275461}

\bibitem[BMS18]{BMS1}
Bhargav Bhatt, Matthew Morrow, and Peter Scholze, \emph{Integral {$p$}-adic
  {H}odge theory}, Publ. Math. Inst. Hautes \'{E}tudes Sci. \textbf{128}
  (2018), 219--397. \MR{3905467}

\bibitem[BMS19]{BMS2}
\bysame, \emph{Topological {H}ochschild homology and integral {$p$}-adic
  {H}odge theory}, Publ. Math. Inst. Hautes \'{E}tudes Sci. \textbf{129}
  (2019), 199--310. \MR{3949030}

\bibitem[BS15]{proetale}
Bhargav Bhatt and Peter Scholze, \emph{The pro-\'{e}tale topology for schemes},
  Ast\'{e}risque (2015), no.~369, 99--201. \MR{3379634}

\bibitem[BS19]{BS19}
Bhargav Bhatt and Peter Scholze, \emph{{Prisms and Prismatic Cohomology}},
  2019, arXiv:1905.08229.

\bibitem[DI87]{DI87}
Pierre Deligne and Luc Illusie, \emph{Rel\`evements modulo {$p^2$} et
  d\'{e}composition du complexe de de {R}ham}, Invent. Math. \textbf{89}
  (1987), no.~2, 247--270. \MR{894379}

\bibitem[Dri18]{drinfeld2018stacky}
Vladimir Drinfeld, \emph{A stacky approach to crystals}, 2018,
  arXiv:1810.11853.

\bibitem[Dri21a]{drinew}
\bysame, \emph{A 1-dimensional formal group over the prismatization of
  {$\mathrm{Spf}~\mathbb{Z}_p$}}, 2021, arXiv:2107.11466.

\bibitem[Dri21b]{prismatization}
\bysame, \emph{Prismatization}, 2021, arXiv:2005.04746.

\bibitem[FM87]{FM87}
Jean-Marc Fontaine and William Messing, \emph{{$p$}-adic periods and {$p$}-adic
  {\'e}tale cohomology}, Current trends in arithmetical algebraic geometry
  ({A}rcata, {C}alif., 1985), Contemp. Math., vol.~67, Amer. Math. Soc.,
  Providence, RI, 1987, pp.~179--207. \MR{902593 (89g:14009)}

\bibitem[Ill71]{Ill}
Luc Illusie, \emph{Complexe cotangent et d\'eformations {I}}, Lecture Notes in
  Mathematics, vol. 239, Springer-Verlag, Berlin, 1971.

\bibitem[Kat87]{Ka87}
Kazuya Kato, \emph{On {$p$}-adic vanishing cycles (application of ideas of
  {F}ontaine-{M}essing)}, Algebraic geometry, {S}endai, 1985, Adv. Stud. Pure
  Math., vol.~10, North-Holland, Amsterdam, 1987, pp.~207--251. \MR{946241}

\bibitem[KP21a]{KP19}
Dmitry Kubrak and Artem Prikhodko, \emph{{Hodge-to-de Rham} degeneration for
  stacks}, 2021, arXiv:1910.12665.

\bibitem[KP21b]{KP21}
\bysame, \emph{$p$-adic {Hodge theory for Artin} stacks}, 2021,
  arXiv:2105.05319.

\bibitem[LL21]{LL20}
Shizhang Li and Tong Liu, \emph{Comparison of prismatic cohomology and derived
  de {Rham} cohomology}, 2021, arXiv:2012.14064.

\bibitem[Lur04]{Lur2}
Jacob Lurie, \emph{Derived algebraic geometry},
  \url{http://people.math.harvard.edu/~lurie/papers/DAG.pdf}, 2004, Ph.d.
  Thesis.

\bibitem[Lur09]{Lu}
Jacob Lurie, \emph{Higher topos theory}, Annals of Mathematics Studies, vol.
  170, Princeton University Press, Princeton, NJ, 2009.

\bibitem[Mon21]{Mon21}
Shubhodip Mondal, \emph{{$\mathbb{G}_a^{\mathrm{perf}}$}-modules and de {Rham}
  cohomology}, 2021, arXiv:2101.03146.

\bibitem[Mon22]{Mon22}
\bysame, \emph{Reconstruction of the stacky approach to de {Rham} cohomology},
  2022, arXiv:2202.07089.

\bibitem[MRT21]{MRT19}
Tasos Moulinos, Marco Robalo, and Bertrand {To{\"e}n}, \emph{{A Universal HKR
  Theorem}}, 2021, arXiv:1906.00118.

\bibitem[Sim96]{Si96}
Carlos Simpson, \emph{Homotopy over the complex numbers and generalized de
  {R}ham cohomology}, Moduli of vector bundles ({S}anda, 1994; {K}yoto, 1994),
  Lecture Notes in Pure and Appl. Math., vol. 179, Dekker, New York, 1996,
  pp.~229--263. \MR{1397992}

\bibitem[{Sta}21]{stacks-project}
The {Stacks project authors}, \emph{The {Stacks} project},
  \url{https://stacks.math.columbia.edu}, 2021.

\bibitem[To{\"{e}}06]{MR2244263}
Bertrand To{\"{e}}n, \emph{Champs affines}, Selecta Math. (N.S.) \textbf{12}
  (2006), no.~1, 39--135. \MR{2244263}

\end{thebibliography}

\end{document}